\documentclass[11pt,letterpaper]{amsart}
\usepackage{url,amsmath,amsthm, amssymb}
\usepackage[margin=1in]{geometry}
\usepackage[pdfencoding=auto, psdextra, colorlinks=true, linkcolor = blue, urlcolor=blue, citecolor=blue, bookmarksdepth=10]{hyperref}
\usepackage{bookmark}
\usepackage[all]{mathtools}
\usepackage{mathrsfs}
\usepackage{mathdots}
\usepackage{amsxtra,amscd}
\usepackage{amsrefs}
\usepackage{enumerate, enumitem}
\usepackage{arydshln}
\usepackage{multirow}
\usepackage{bm} 
\usepackage[all,arc]{xy}

 \usepackage{amssymb,amsrefs}
\usepackage{amsmath}
\usepackage{amscd,latexsym}
\usepackage{bookmark}
 \usepackage{amssymb,amsfonts,bm}
\usepackage[all,arc]{xy}
\usepackage{enumerate}
\usepackage{mathrsfs}
\usepackage{amscd}
\usepackage{tikz}
\usepackage{tikz-cd}
\usepackage{amsmath}
\usepackage{latexsym}
\usepackage{mathdots}
\usepackage{amsrefs}
\usepackage{graphicx}
\usepackage{mathtools}
\usepackage{leftidx}
\usepackage{tensor}
\usepackage{dsfont}
\usepackage{tikz-cd} 
\usepackage[new]{old-arrows}
\usepackage{xcolor}
\usepackage{bbm}
\usepackage{enumitem}

\allowdisplaybreaks

\theoremstyle{plain}
\newtheorem{theorem}{Theorem}[section]
\newtheorem{corollary}[theorem]{Corollary}

\newtheorem{proposition}[theorem]{Proposition}
\numberwithin{equation}{section}
\newtheorem*{whypothesis}{Working Hypothesis}

\theoremstyle{definition}
 \newtheorem{definition}[theorem]{Definition}

\theoremstyle{remark}

\newtheorem*{acknowledgements}{Acknowledgments}

\begin{document}

\title[The local exterior square and Asai $L$-functions]{The local exterior square and Asai $L$-functions for $GL(n)$ in odd characteristic}

\author[]{Yeongseong Jo}
\address{Department of Mathematics Education, Ewha Womans University, Seoul 03760, Republic of Korea}
\email{\href{mailto:jo.59@buckeyemail.osu.edu}{jo.59@buckeyemail.osu.edu}}


\begin{abstract}
Let $F$ be a non-archimedean local field of odd characteristic $p > 0$. In this paper, we consider
local exterior square $L$-functions $L(s,\pi,\wedge^2)$, Bump-Friedberg $L$-functions $L(s,\pi,BF)$, and Asai $L$-functions $L(s,\pi,As)$
of an irreducible admissible representation $\pi$ of $GL_m(F)$. In particular, we establish that those $L$-functions, 
via the theory of integral representations, are equal to their corresponding Artin $L$-functions 
  $L(s,\wedge^2(\phi(\pi)))$, $L(s+1/2,\phi(\pi))L(2s,\wedge^2(\phi(\pi)))$, and  $L(s,As(\phi(\pi)))$ of the associated Langlands parameter $\phi(\pi)$
under the local Langlands correspondence. These are achieved by proving the identity for irreducible supercuspidal representations,
exploiting the local to global argument due to Henniart and Lomeli. 
 \end{abstract}

\maketitle

\section{Introduction}

Let $F$ be a non-archimedean local field of positive characteristic $p \neq 2$. Let $\pi$ be an irreducible 
admissible representation of $GL_m(F)$, where $m$ is a positive integer. 
The local Langlands correspondence provides a bijection between the set of 
equivalence classes of irreducible admissible (complex-valued) representations of $GL_m(F)$ 
and the set of equivalence classes of $m$-dimensional Weil-Deligne representations of Weil-Deligne group $W'_F$ of $F$.
Let $r$ denote either the exterior square representation $\wedge^2 :  GL_m(\mathbb{C}) \rightarrow  GL_{m(m-1)/2}(\mathbb{C})$ or the Asai representation (the twisted tensor induction) $As : GL_m(\mathbb{C}) \rightarrow GL_{m^2}(\mathbb{C})$ (Cf. \cite[\S 2.1]{AR}, \cite[\S 1.2]{Shankman}) of 
the Langlands dual group $GL_m(\mathbb{C})$ of $GL_m(F)$. Let $L(s,(r \circ \phi)(\pi) ), s \in \mathbb{C}$ be the 
Artin local $L$-function defined by Deligne and Langlands \cite{Tate}, where $\phi(\pi) : W'_F \rightarrow GL_m(\mathbb{C})$ is the Weil-Deligne representation
corresponding $\pi$ under the local Langlands correspondence. In this paper, we address that in the cases of
exterior square, Bump-Friedberg, and Asai local $L$-functions, $L$-factors are compatible with the local Langlands correspondence, and establish a series of equalities of local $L$-functions:
\[
\begin{aligned}
&\bullet  \text{(Jacquet-Shalika cases, Theorem \ref{exter-equal})} & &L(s,\pi,\wedge^2)&&=&&L(s,\wedge^2(\phi(\pi)));\\
& \bullet \text{(Bump-Friedberg cases, Theorem \ref{BF-equal})} &  &L(s,\pi,BF)&&=&&L(s+1/2,\phi(\pi))L(2s,\wedge^2(\phi(\pi)));\\
 &\bullet \text{(Flicker cases, Theorem \ref{equal-asai})} &  &L(s,\pi,As)&&=&&L(s,As(\phi(\pi)),\\
\end{aligned}
\]
where $L$-factors on the left hand sides are defined by the theory of integral representations in positive characteristic. 

\par
In late 1980's, global Zeta integrals for  (partial) Asai $L$-functions for $GL_m$ appeared in the paper by Flicker \cite{Fil88,Fil93}. Around that time, 
Jacquet and Shalika \cite{JaSh88} and Bump and Friedberg \cite{BF} independently constructed two different integral representations for an (incomplete) exterior square $L$-function associated to a cuspidal automorphic representation on $GL_m$ over a global field. Recently there has been renowned interest in the local theory of Asai and exterior square
$L$-functions via Rankin-Selberg methods. In characteristic $0$, the identities were shown for Jacquet-Shalika integrals by the author \cite{Jo20-2}, and Bump-Friedberg integrals \cite{Matringe15} and Flicker integrals \cite{Matringe09,Matringe11} by Matringe. As a matter of fact, these results improve discrete series cases of Kewat and Raghunathan \cite{KeRa} for Jacquet-Shalika integrals and Anandavardhanan and Rajan \cite{AR} for Flicker integrals. In the positive characteristic $p > 0$, Artin $L$-factors $L(s,\wedge^2(\phi(\pi)))$ and $L(s,As(\phi(\pi))$ coincide with the corresponding ones $\mathcal{L}(s,\pi,\wedge^2)$ and $\mathcal{L}(s,\pi,As)$ via the Langlands-Shahidi method by a sequence of work of Henniart and Lomel\'{\i} \cite{HL-Exterior,HL-Asai,HL-RS}. Similar problems have been worked out by Henniart \cite{Henniart} in the characteristic zero cases.

\par
The method to prove the matching is developed by Cogdell and Piatetski-Shapiro \cite{CogPS} in the framework of local $L$-functions of pairs of irreducible generic representations $(\pi_1,\pi_2)$. The computation of local Rankin-Selberg $L$-functions is boiled down to decomposing it as the product of what is called the exceptional $L$-functions (in the sense of \cite{CogPS}) $L_{ex}(s,\pi^{(k_1)}_1 \times \pi_2^{(k_2)})$ for pairs of Bernstein-Zelevinsky's derivatives $(\pi_1^{(k_1)},\pi_2^{(k_2)})$. The ``derivative" in the sense of Bernstein-Zelevinsky $\pi^{(k)}$ is given by representations of smaller groups $GL_{m\text{-}k}(F)$. The advantage of adapting such derivatives enables us to
proceed it by induction on the rank $m$-$k$ of general linear groups $GL_{m\text{-}k}(F)$.

\par
Each pole of exceptional $L$-functions, which we often refer to ``exceptional poles", is astonishingly characterized by local distinctness or existence of certain models. Sooner and later, the classification of irreducible generic distinguished representations attains its own merit and it is widely explored in various work. Indeed, topics of the classification are brought to light
for $(S_{2n}(F),\Theta)$-distinguished representations (Shalika models) in \cite{Matringe17}, for $H_m(F)$-distinguished representation (Linear / Friedberg-Jacquet models)  in \cite{Matringe15},
for $GL_m(F)$-distinguished representations (Flicker-Rallis models) in \cite{Matringe11} by Matringe, and for $(GL_m(F), \theta)$-distinguished representations (Bump-Ginzburg models) in \cite{Kaplan} by Kaplan. In particular, the part of Theorems is highlighted as, so to speak, ``the hereditary property of models" motivated from the classification of irreducible admissible generic representations of $GL_m(F)$ (Whittaker models) due to Rodier \cite[P. 427]{Rodier}.

\par
When we combine these two ingredients, the factorization of local $L$-functions and the classification of distinguished representations, we are guided to obtain major applications, the inductive relation of local $L$-functions for irreducible generic representations $\pi=\mathrm{Ind}^{GL_m}_{\rm Q}(\Delta_1 \otimes \Delta_2 \otimes \dotsm \otimes \Delta_t)$ (Corollary \ref{generic}):
\begin{equation}
\label{L-equal}
  L(s,\pi,\wedge^2)=\prod_{1 \leq k \leq t} L(s,\Delta_k,\wedge^2) \prod_{1 \leq i < j  \leq t} L(s, \Delta_i \times \Delta_j),
\end{equation}
and the weak multiplicativity of $\gamma$-factors for parabolically normalized induced 
(not necessarily irreducible) representations $\pi=\mathrm{Ind}^{GL_m}_{\rm Q}(\Delta_1 \otimes \Delta_2 \otimes \dotsm \otimes \Delta_t)$ (Theorem \ref{Deform-gamma}):
\begin{equation}
\label{gamma-equal}
\gamma(s,\pi,\wedge^2,\psi)
\sim \prod_{1 \leq k \leq t} \gamma(s,\Delta_k,\wedge^2,\psi) \prod_{1 \leq i < j  \leq t} \gamma(s, \Delta_i \times \Delta_j,\psi),
\end{equation}
where $\sim$ means the equality up to a unit in $\mathbb{C}[q^{\pm s}]$ and $\Delta_i$'s are discrete series representations. 
Building upon \eqref{L-equal}, we can incorporate the Langlands classification of irreducible admissible representations in terms of discrete series ones into the theory of local $L$-functions. In turn, \eqref{gamma-equal} allows us to compute local $L$-functions further in accordance with the Bernstein-Zelevinsky classification of discrete series representations in terms of irreducible supercuspidal ones. As a consequence, we express all exterior square $L$-factors for irreducible admissible representations in terms of $L$-factors for irreducible supercuspidal ones in a purely local mean. This comes down to reducing our main questions to all irreducible supercuspidal representations, which are eventually served as building blocks.

\par
We emphasize that the third identity in Flicker cases is not new and it can be found in the Appendix A of \cite{AKMSS}. However, we discovered that
our technique seems (at least to us) to carry out uniformly to other $L$-function for $GL_m$. As an application of our approach,
the main result of \cite{AKMSS} immediately implies the agreement of local Asai $L$-factors for irreducible supercuspidal representations, which is sufficient
to extend it to all irreducible admissible representations, reflecting on the local Langlands correspondence. At this point, unlike \cite[Appendix A]{AKMSS}, additional globalizations are not required to
generalize the equality unconditionally.  In the course of following the direction taken by \cite[Appendix A]{AKMSS}, we encounter a few stumbling blocks.
In contrast to characteristic 0 cases, we could not find a good way of adjusting the globalization of discrete series representations in Gan-Lomel\'{\i} \cite[Proposition 8.2]{Gan-L} 
to our circumstance. As seen in several others work \cite{AR,KeRa,Kable,Yamana}, there might be no guarantee that the different places $v_1$ and $v_2$ are co-prime in order to conclude that $\log(q_{v_1})/\log(q_{v_2})$ is irrational, and this coprimity condition may prompt an issue in characteristic $p > 0$. We propose to resolve all these difficulties by
globalizing irreducible supercuspidal representations in Henniart-Lomel\'{\i}  \citelist{\cite{HL-Exterior}\cite{HL-RS}*{Theorem 3.1}} and controlling all but one place in which we are interested.

\par
In practice, we 
demonstrate the identity sequently for irreducible supercuspidal representations and eventually for discrete series representations under the Working hypothesis analogous to Kaplan's enquiry \cite[Remark 4.18]{Kaplan} that $(GL_m(F), \theta)$-distinguished discrete series representations in positive characteristic are self-dual. Thankfully, we remove the hypothesis by investigating irreducible generic subquotients of principal series representations. We expect to overcome Kaplan's issue beyond the principal series cases by reconciling the different definitions of local symmetric square $L$-functions possessing their own insights about representations. The poles of $L(s,\pi,\mathrm{Sym}^2)$  by means of the Rankin-Selberg method can be determined by the occurrence of $(GL_m(F), \theta)$-distinguished representations \cite{Yamana}, whereas the symmetric square $L$-functions through the Langlands-Shahidi method $\mathcal{L}(s,\pi,\mathrm{Sym}^2)$ can be related to the presence of the self-duality $\pi \simeq \widetilde{\pi}$, using the Rankin-Selberg $L$-factor $L(s,\pi \times \pi)$ as a product of $\mathcal{L}(s,\pi,\wedge^2)$ and $\mathcal{L}(s,\pi,\mathrm{Sym}^2)$ \cite{HL-Exterior,HL-RS}. Taking it for granted that  $L(s,\pi,\mathrm{Sym}^2)$ can be factored in terms of exceptional $L$-factors for derivatives (Cf. \cite[Theorem 3.15]{Jo21}), our discourse sheds light on some impetus toward systematic development of symmetric square $L$-factors via integral representations \cite{Yamana} in number theoretic aspects and the classification of $(GL_m(F), \theta)$-distinguished representations over local function fields \cite{Kaplan} in representation theoretic perspectives. We will return these matters in near future.

\par
Finally, it is worth pointing out that the result of this paper will be used to supply the proof of the clam in the preprint by Chen and Gan \cite[Theorem 1.1]{Chen-Gan} that {\it the exterior square $L$-function 
can be equivalently defined by the Langlands-Shahidi method or the local zeta integrals of Jacquet-Shalika} \cite{JaSh88} in positive characteristic.

\par
Let us overview the content of this paper. Section \ref{sec2.1} begins with a summary of the theory of derivatives of Bernstein and Zelevinsky and the basic existence Theorem of Jacquet-Shalika integrals. Section \ref{sec:2.2} is devoted to classifying all irreducible generic distinguished representation with respect to given closed algebraic subgroups, especially $H_{2n}$ and $S_{2n}$, due to Matringe. 
By combining the factorization of Section \ref{sec2.1}, the classification of Section \ref{sec:2.2}, and the method of deformations and specializations, we prove a weak version of multiplicativity of $\gamma$-factors and the inductive relation of $L$-factors. Using the globalization of irreducible supercuspidal representation presented in Section \ref{sec3.2}, if necessary, we complete computing local exterior square $L$-functions at the end of Section \ref{sec3.2}, local Bump-Friedberg $L$-functions in Section \ref{sec:BF}, and local Asai $L$-functions in Section \ref{sec:Asai}.

\section{Jacquet-Shalika Zeta Integrals}
\label{sec2}

\subsection{Derivatives and exceptional poles}
\label{sec2.1}
Let $F$ be a non-archimedean local field of characteristic $p \neq 0,2$. We let $\mathcal{O}$ denote its ring of integers, $\mathfrak{p}$ its maximal ideal, $q$ the cardinality of its residual field. We will let $\varpi$ denote a uniformizer, so $\mathfrak{p}=(\varpi)$. We normalize the absolute value by $|\varpi|^{-1}=|\mathcal{O} \slash \mathfrak{p}|$. The character of $GL_m$ given by $g \mapsto |\mathrm{det}(g)|$ is denoted by $\nu$.

\par
 For the group $GL_m:=GL_m(F)$, we often confront the two cases: $m$ is even and $m$ is odd. For the former we let $m=2n$, and for the latter $m=2n+1$. Let $\sigma_m$ be the permutation matrix given by
\[
 \sigma_{2n}=\begin{pmatrix} 1 & 2 & \dotsm & n & | & n+1 & n+2 & \dotsm & 2n \\ 
                                                1 & 3 & \dotsm & 2n-1 & | &  2 & 4 & \dotsm &2n \\ 
                                                \end{pmatrix}
\]
when $m=2n$ is even, and by
\[
 \sigma_{2n+1}=\begin{pmatrix} 1 & 2 & \dotsm & n & | & n+1 & n+2 & \dotsm & 2n & 2n+1 \\ 
                                                1 & 3 & \dotsm & 2n-1 & | &  2 & 4 & \dotsm &2n & 2n+1 \\ 
                                                \end{pmatrix}
\]
when $m=2n+1$ is odd. Let $B_m$ be the Borel subgroup consisting of the upper triangular matrices with Levi subgroup $A_m$ of diagonal matrices and unipotent radical $N_m$. We write $Z_m$ denote to the center consisting of scalar matrices. We define $P_m$ the mirabolic subgroup given by
\[
 P_m=\left\{ \begin{pmatrix} g & {^tu} \\ & 1 \end{pmatrix} \;\middle| \; g \in GL_{m-1}, u \in F^{m-1}  \right\}.
\]
We denote by $U_m$ the unipotent radical of $P_m$. As a group, $P_m$ has a structure of a semi-direct product $P_m=GL_{m-1} \ltimes U_m$.
We let $\mathcal{M}_m$ be the $m \times m$ matrices, $\mathcal{N}_m$ the subspace of upper triangular matrices of $\mathcal{M}_m$. Let $\{ e_i \,|\, 1 \leq i \leq m\}$ be the standard low basis of $F^m$.

\par
We let $\psi_F$ denote a non-trivial additive character of $F$. We let $\psi$ denote the character of $N_m$ defined by
\[
  \psi(n)=\psi_F ( \sum_{i=1}^{n-1} n_{i,i+1} ), \quad n=(n_{i,j}) \in N_m.
\]
We denote by $\mathcal{A}_F(m)$ the set of equivalence classes of all admissible representations of $GL_m$ on complex vector spaces. 
Furthermore we say that a representation $\pi \in \mathcal{A}_F(m)$ is called {\it generic} if $\mathrm{Hom}_{N_m}(\pi,\psi) \neq \{ 0 \}$.
We say that a representation $\pi \in \mathcal{A}_F(m)$ is of {\it Whittaker type} if
\[
  \mathrm{dim}_{\mathbb{C}} \mathrm{Hom}_{N_m}(\pi,\psi)=1. 
\]
For any character $\chi$ of $F^{\times}$, $\chi$ can be uniquely decomposed as $\chi=\chi_0\nu^{s_0}$, where $\chi_0$ is a unitary character and $s_0$ is a real number. We use the notation $s_0=\mathrm{Re}(\chi)$ for the real part of the exponent of the character $\chi$.

\par
If $\pi \in \mathcal{A}_F(m)$ is irreducible and generic, it is known from \cite{GK} that $\pi$ is of Whittaker type. By Frobenius reciprocity, there exists a unique embedding of $\pi$ into $\mathrm{Ind}^{GL_m}_{N_m}(\psi)$ up to scalar. The image $\mathcal{W}(\pi,\psi)$ of $V_{\pi}$ is called the {\it Whittaker model} of $\pi$. For a non-zero functional $\lambda \in \mathrm{Hom}_{N_m}(\pi,\psi)$,
we define  the {\it Whittaker function} $W_v \in \mathcal{W}(\pi,\psi)$ associated to $v \in V_{\pi}$ by 
\[
  W_v(g)=\lambda(\pi(g)v), \quad g \in GL_m.
\]
We set $W:=W_v$. It follows from \cite[Lemma 4.5]{BeZeSurvey} and \cite[\S 9]{Zelevinsky} that if $\Delta_1,\Delta_2,\dotsm,\Delta_t$ are irreducible {\it essentially square integrable}, which we call {\it discrete series representations}, then
the representation of the form $\mathrm{Ind}^{GL_m}_{\rm Q}(\Delta_1 \otimes \Delta_2 \otimes \dotsm \otimes \Delta_t)$ is a representation of Whittaker type, where the induction is the normalized parabolic induction from the standard parabolic subgroup $\rm Q$ attached to the partition $(m_1,m_2,\dotsm,m_t)$ of $m$ and $\Delta_i \in \mathcal{A}_F(m_i)$. Also, whenever the parabolic subgroup ${\rm Q}$ and ambient group $GL_m$ are clear from the context, we simply write $\mathrm{Ind}(\Delta_1 \otimes \Delta_2 \otimes \dotsm \otimes \Delta_t)$.

\par
Let $\mathrm{Rep}(G)$ denote the category of smooth representations of an $l$-group $G$. There are four functors $\Psi^-$, $\Psi^+$, $\Phi^-$, and $\Phi^+$. The functor $\Psi^-$ is
a normalized Jacquet functor and  $\Phi^-$ a normalized $\psi$-twisted Jacquet functor from $\mathrm{Rep}(P_m)$ to $\mathrm{Rep}(GL_{m-1})$ and $\mathrm{Rep}(P_{m-1})$, respectively. Given $\tau \in \mathrm{Rep}(P_m)$ on the space $V_{\tau}$, $\Psi^-(\tau)$ is realized on the space $V_{\tau} \slash V_{\tau}(U_m, {\bf 1})$ with the action 
$\Psi^-(\tau)(g)(v+V_{\tau}(U_m, {\bf 1}))=|\mathrm{det}(g)|^{-\frac{1}{2}}(\tau(g)v+V_{\tau}(U_m, {\bf 1}))$ and the subspace $V_{\tau}(U_m, {\bf 1}))=\langle \tau(u)v-v \; | \; v \in V_{\tau}, u \in U_m \rangle$. Likewise $\Phi^-(\tau)$ is realized on the space $V_{\tau} \slash V_{\tau}(U_m, \psi)$ with the action 
$\Phi^-(\tau)(p)(v+V_{\tau}(U_m, \psi))=|\mathrm{det}(p)|^{-\frac{1}{2}}(\tau(p)v+V_{\tau}(U_m, \psi))$ and the subspace $V_{\tau}(U_m, \psi))=\langle \tau(u)v-\psi(u)v \; | \; v \in V_{\tau}, u \in U_m \rangle$. The functor  $\Psi^+$ and  $\Phi^+$ are normalized compactly supported inductions from $\mathrm{Rep}(GL_{m-1})$ and $\mathrm{Rep}(P_{m-1})$, respectively to $\mathrm{Rep}(P_m)$. Given $\sigma \in \mathrm{Rep}(GL_{m-1})$, $ \Psi^+(\sigma)=\mathrm{ind}^{P_m}_{GL_{m-1}U_m}(|\mathrm{det}(g)|^{\frac{1}{2}}\sigma \otimes {\bf 1} )=|\mathrm{det}(g)|^{\frac{1}{2}}\sigma \otimes {\bf 1}$ is realized on the space $V_{\sigma}$, where ``ind"  is a compactly supported induction. If $\sigma \in \mathrm{Rep}(P_{m-1})$, then
$ \Phi^+(\sigma)=\mathrm{ind}^{P_m}_{P_{m-1}U_m}(|\mathrm{det}(g)|^{\frac{1}{2}}\sigma \otimes \psi )$.

\par
For $\tau \in \mathrm{Rep}(P_m)$, four functors are utilized to define what is called {\it Bernstein-Zelevinsky $k^{th}$-derivatives} $\tau^{(k)}$. Let $\tau^{(k)} \in \mathrm{Rep}(GL_{m-k})$ be $\tau^{(k)}=\Psi^-(\Phi^-)^{k-1}(\tau)$ for $1 \leq k \leq m$. The smooth representation $\tau$ affords a natural filtration by $P_m$-modules 
\[
 0 \subseteq \tau_m \subseteq \tau_{m-1} \subseteq \dotsm \subseteq \tau_1=\tau
\]
such that $\tau_k \slash \tau_{k+1}=(\Phi^+)^{k-1} \Psi^+(\tau^{(k)})$ and $\tau_k=(\Phi^+)^{k-1}(\Phi^-)^{k-1}(\tau)$. Let $\pi=\mathrm{Ind}^{GL_m}_{\rm Q}(\Delta_1 \otimes \dotsm \otimes \Delta_t)$ be a parabolically induced representation, where $\Delta_i$ is an irreducible essentially square integrable representation of $GL_{m_i}$ so that $m=m_1+\dotsm+m_t$.
Then $\pi^{(k)}$ has a filtration whose successive quotients are isomorphic to $\mathrm{Ind}(\Delta_1^{(k_1)} \otimes \dotsm \otimes \Delta_t^{(k_t)})$ with $k=k_1+\dotsm+k_t$ \cite[Theorem 4.4, Lemma 4.5]{BeZe}. For every $0 \leq k \leq m-1$, let $(\omega_{\pi_{i_k}^{(k)}})_{i_k=1,2,\dotsm,r_k}$ be the family of the central characters of non-zero successive quotient of the form $\pi_{i_k}^{(k)}=\mathrm{Ind}(\Delta_1^{(k_1)} \otimes \dotsm \otimes \Delta_t^{(k_t)})$.

\par
Let $\mathcal{S}(F^n)$ be the space of smooth locally constant compactly supported functions on $F^n$. For each Whittaker function $W \in \mathcal{W}(\pi,\psi)$ and Schwartz-Bruhat function $\Phi \in \mathcal{S}(F^n)$, we define Jacquet-Shalika integrals:
\[
J(s,W,\Phi)=
 \int_{N_n \backslash GL_n} \int_{\mathcal{N}_n \backslash \mathcal{M}_n} W \left( \sigma_{2n} \begin{pmatrix} I_n & X \\ & I_n  \end{pmatrix} \begin{pmatrix} g & \\ & g \end{pmatrix} \right)  \psi^{-1}({\rm Tr}(X)) \Phi(e_ng) |\mathrm{det}(g)|^s dX dg
\]
in the even case $m=2n$ and 
\[
\begin{split}
J(s,W,\Phi)=& \int_{N_n \backslash GL_n} \int_{\mathcal{N}_n \backslash \mathcal{M}_n} \int_{F^n} W \left( \sigma_{2n+1} \begin{pmatrix} I_n &X& \\ &I_n& \\ &&1 \end{pmatrix} \begin{pmatrix} g && \\ &g& \\ &&1 \end{pmatrix} \begin{pmatrix} I_n && \\ &I_n& \\ &y&1 \end{pmatrix} \right) \\
 &\times \psi^{-1}({\rm Tr}(X)) \Phi(y) |\mathrm{det}(g)|^{s-1}  dy dX dg\\
\end{split}
\]
in the odd case $m=2n+1$. Several nice consequences follow from inserting asymptotic formula over the torus for $W$ into the local Zeta integral $J(s,W,\Phi)$ \cite[Theorem 3.3, Lemma 3.10]{Jo20}.

\begin{theorem}
\label{exterior-structure}
Let $\pi=\mathrm{Ind}^{GL_m}_{\rm Q}(\Delta_1 \otimes \Delta_2 \otimes \dotsm \otimes \Delta_t)$ be a parabolically induced representation. 
Let $W \in \mathcal{W}(\pi,\psi)$ and $\Phi \in \mathcal{S}(F^n)$.   
\begin{enumerate}
\item[$(\mathrm{i})$]\hspace{-2mm}-$(1)$ $($Even cases $m=2n$$)$ If we have,
$
  \mathrm{Re}(s) > -\frac{1}{k} \omega_{\pi^{(2n-2k)}_{i_{2n-2k}}},
$
for all $1 \leq k \leq n$ and all $1 \leq i_{2k} \leq r_{2k}$, then each local integral $J(s,W,\Phi)$ converges absolutely.
\item[$(\mathrm{i})$]\hspace{-2mm}-$(2)$ $($Odd cases $m=2n+1$$)$ If we have
$
 \mathrm{Re}(s) > -\frac{1}{k} \omega_{\pi^{(2n+1-2k)}_{i_{2n+1-2k}}},
$
for all $1 \leq k \leq n$ and all $1 \leq i_{2k-1} \leq r_{2k-1}$, then each local integral $J(s,W,\Phi)$ converges absolutely.
\item[$(\mathrm{ii})$] Each $J(s,W,\Phi)$ is a rational function in $\mathbb{C}(q^{-s})$, hence $J(s,W,\Phi)$ as a function of $s$ extends meromorphically to all $\mathbb{C}$. 
\item[$(\mathrm{iii})$] Each $J(s,W,\Phi)$ can be written with a common denominator determined by $\pi$. Hence the family has ``bounded denominators".
\end{enumerate}
\end{theorem}

Let $\mathcal{J}(\pi)$ be the complex linear space of the local integrals $J(s,W,\Phi)$. The family of local integrals $\mathcal{J}(\pi)$ is a $\mathbb{C}[q^{\pm s}]$-fractional ideal of $\mathbb{C}(q^{-s})$ containing $1$ \cite[Theorem 3.6, Theorem 3.9]{Jo20}. Since the ring $\mathbb{C}[q^s,q^{-s}]$ is a principal ideal domain, the fractional ideal $\mathcal{J}(\pi)$ has a generator. Since $1 \in \mathcal{J}(\pi)$ we can take a generator having numerator $1$ and normalized (up tp units) to be of the form $P(q^{-s})^{-1}$ with $P(X) \in \mathbb{C}[X]$ having $P(0)=1$. The {\it local exterior square $L$-function} or simply {\it exterior square $L$-factor}
\[
  L(s,\pi,\wedge^2)=\frac{1}{P(q^{-s})}
\]
is defined to be the normalized generator of the fractional ideal $\mathcal{J}(\pi)$ spanned by the local Zeta integrals $J(s,W,\Phi)$.
 
\par 
 We define the Fourier transform on $\mathcal{S}(F^m)$ by
 \[
   \hat{\Phi}(y)=\int_{F^n} \Phi(x) \psi (x\,{^ty})\,dx.
 \]
 We assume that the measure on $F^m$ is the self-dual measure. Then the Fourier inversion takes the form $\hat{\hat{\Phi}}(x)=\Phi(-x)$. Let 
 \[
  w_m:=\begin{pmatrix}  \vspace{-2ex} && 1 \\ & \iddots& \vspace{-2ex} \\ 1 && \end{pmatrix}
 \]
 denote the long Weyl element in $GL_m$. For $(\pi,V_{\pi}) \in \mathrm{Rep}(GL_m)$, let $\pi^{\iota}$ denote the representation of $GL_m$ on the same space $V_{\pi}$ given by $\pi^{\iota}(g)=\pi(^tg^{-1})$. If $\pi$ is irreducible, it is known that $\pi^{\iota} \simeq \widetilde{\pi}$, the contragredient representation of $\pi$. The parabolically induced representation $\pi^{\iota}=\mathrm{Ind}(\widetilde{\Delta}_t \otimes \widetilde{\Delta}_{t-1} \otimes \dotsm \otimes \widetilde{\Delta}_{1}  )$ is again of Whittaker type. If $W \in \mathcal{W}(\pi,\psi)$, then $\widetilde{W}(g):=W(w_m\,{^tg^{-1}})$ belongs to $\mathcal{W}(\pi^{\iota},\psi^{-1})$. We let $\tau_m$ be a matrix given by $\begin{pmatrix} & I_n \\ I_n & \end{pmatrix}$ when $m=2n$ is even, and by $ \begin{pmatrix} & I_n &\\ I_n && \\ &&1 \end{pmatrix}$ when $m=2n+1$ is odd. As a consequence of the uniqueness of bilinear forms on $\mathcal{W}(\pi,\psi) \times \mathcal{S}(F^n)$, we can define the local $\gamma$-factor which gives rise to the local functional equation for our integrals $J(s,W,\Phi)$ \citelist{\cite{CoMa}\cite{Matringe}} (Cf. \cite[Theorem 2.10, (2.1)]{Jo20-2}).

\begin{theorem}
\label{exterior-func}
 Let $\pi=\mathrm{Ind}^{GL_m}_{\rm Q}(\Delta_1 \otimes \Delta_2 \otimes \dotsm \otimes \Delta_t)$ be a parabolically induced representation of $GL_m$. 
Then there is a rational function $\gamma(s,\pi,\wedge^2,\psi) \in \mathbb{C}(q^{-s})$  such that for every $W$ in $\mathcal{W}(\pi,\psi)$, and every $\Phi$ in $\mathcal{S}(F^n)$, we have 
\[
 J(1-s,\varrho(\tau_m)\widetilde{W},\hat{\Phi})=\gamma(s,\pi,\wedge^2,\psi)J(s,W,\Phi),
\]
where $\varrho$ denotes right translation.
\end{theorem}

An equally important local factor is the local $\varepsilon$-factor
\[
  \varepsilon(s,\pi,\wedge^2,\psi)=\gamma(s,\pi,\wedge^2,\psi) \frac{L(s,\pi,\wedge^2)}{L(1-s,\pi^{\iota},\wedge^2)}
\]
which is an invertible element $\varepsilon(s,\pi,\wedge^2,\psi)$ in $\mathbb{C}[q^{\pm s}]$. With the local $\varepsilon$-factor the functional equation becomes
\[
  \frac{J(1-s,\varrho(\tau_m)\widetilde{W},\hat{\Phi})}{L(1-s,\pi^{\iota},\wedge^2)}
  =\varepsilon(s,\pi,\wedge^2,\psi) \frac{J(s,W,\Phi)}{L(s,\pi,\wedge^2)}.
\]

\par
 Let $\pi=\mathrm{Ind}^{GL_{2n}}_{\rm Q}(\Delta_1 \otimes \Delta_2 \otimes \dotsm \otimes \Delta_t)$ be a parabolically induced representation. 
 Let $\mathcal{S}_0(F^n)$ be the subspace of $\Phi \in \mathcal{S}(F^n)$ for which $\Phi(0,0,\dotsm,0)=0$. Suppose that there exists a function in $\mathcal{J}(\pi)$ having a pole of order $d_{s_0}$ at $s=s_0$. We investigate the rational function defined by an individual Zeta integral $J(s,W,\Phi)$. Then the Laurent expansion about $s=s_0$ will take the form
 \[
  J(s,W,\Phi)=\frac{B_{s_0}(W,\Phi)}{(q^s-q^{s_0})^{d_{s_0}}}+ \text{higher order terms}.
 \]
 We define the Shalika subgroup $S_{2n}$ of $GL_{2n}$ by
 \[
   S_{2n}=\left\{ \begin{pmatrix} I_n & Z \\ & I_n \end{pmatrix} \begin{pmatrix} h & \\ & h \end{pmatrix} \; \middle| \; Z \in \mathcal{M}_n,\; h \in GL_n  \right\}.
 \]
Let us denote an action of the Shalika subgroup $S_{2n}$ on $\mathcal{S}(F^n)$ by
\[
  R \left(  \begin{pmatrix} I_n & Z \\ & I_n \end{pmatrix} \begin{pmatrix} h & \\ & h \end{pmatrix}  \right) \Phi(x)=\Phi(xh)
\]
 for $\Phi \in \mathcal{S}(F^n)$. The coefficient of the leading term, $B_{s_0}(W,\Phi)$, will define a non-trivial bilinear form on $\mathcal{W}(\pi,\psi) \times \mathcal{S}(F^n)$
 enjoying the quasi-invariance
 \[
  B_{s_0} \left( \varrho (g)W, R (g) \Phi   \right)=|\mathrm{det}(h)|^{-s_0}\psi(\mathrm{Tr}(Z)) B_{s_0}(W,\Phi) \quad \text{for} \quad g=\begin{pmatrix} I_n & Z \\ & I_n \end{pmatrix} \begin{pmatrix} h & \\ & h \end{pmatrix}  \in S_{2n}.
 \]
 The pole at $s=s_0$ of the family $\mathcal{J}(\pi)$ is called {\it exceptional} if the associated bilinear form $B_{s_0}(W,\Phi)$
 vanishes identically on $\mathcal{W}(\pi,\psi) \times \mathcal{S}_0(F^n)$. If $s=s_0$ is an exceptional pole of $\mathcal{J}(\pi)$ then the bilinear form $B_{s_0}$
 factors to a non-zero bilinear form on $\mathcal{W}(\pi,\psi) \times \mathcal{S}(F^n) \slash \mathcal{S}_0(F^n)$. The quotient $\mathcal{S}(F^n) \slash \mathcal{S}_0(F^n)$
 is isomorphic to $\mathbb{C}$ via the map $\Phi \mapsto \Phi(0)$. Let $\Theta$ be the character of $S_{2n}$ given by
 \[
  \Theta \left( \begin{pmatrix} I_n & Z \\ & I_n \end{pmatrix} \begin{pmatrix} h & \\ & h \end{pmatrix}  \right)=\psi(\mathrm{Tr}(Z)).
 \]
We say that $\pi \in \mathcal{A}_F(2n)$ is called $(S_{2n},\Theta)$-{\it distinguished} if ${\rm Hom}_{S_{2n}}(\pi,\Theta) \neq \{ 0\}$. A nonzero linear functional $\Lambda$ in ${\rm Hom}_{S_{2n}}(\pi,\Theta)$ ($\Lambda_s$ in ${\rm Hom}_{S_{2n}}(\pi,\nu^{-\frac{s}{2}}\Theta)$) is called a {\it Shalika functional} (a {\it twisted Shalika functional}, respectively). If $s=s_0$ is an exceptional pole, then the form $B_{s_0}$ can be written as $B_{s_0}(W,\Phi)=\Lambda_{s_0}(W)\Phi(0)$ with $\Lambda_{s_0}$ the Shalika functional on $\mathcal{W}(\pi,\psi)$.
Using the notation, we let
 \[
   L_{ex}(s,\pi,\wedge^2)=\prod_{s_0} (1-q^{s_0}q^{-s})^{d_{s_0}},
 \]
where $s_0$ runs through all the exceptional poles of $\mathcal{J}(\pi)$ with $d_{s_0}$ the maximal order of the pole at $s=s_0$. The factorization of local exterior square $L$-functions proposed by Cogdell and Piatetski-Shapiro asserts that it can be expressed in terms of the exceptional exterior square $L$-factors of the derivatives of $\pi$ \cite{Jo20}.

 \begin{theorem}
 \label{factor}
 Let $\pi=\mathrm{Ind}^{GL_m}_{\rm Q}(\Delta_1 \otimes \Delta_2 \otimes \dotsm \otimes \Delta_t)$ be an irreducible generic representation of $GL_m$ such that 
 all the derivatives $\pi^{(k)}$ of $\pi$ are completely reducible with irreducible generic constituents of the form $\pi^{(k)}_{i}=\mathrm{Ind}(\Delta_1^{(k_1)} \otimes \dotsm \otimes \Delta_t^{(k_t)})$ with $k=k_1+\dotsm+k_t$. For each $k$, $i$ is indexing the partition of $k$. Then
 \begin{enumerate}[label=$(\mathrm{\roman*})$]
\item\label{Even-Facto} $($Even cases $m=2n$$)$ $L(s,\pi,\wedge^2)=\underset{k,\; i}{l.c.m.}\{ L_{ex}(s,\pi^{(2k)}_{i},\wedge^2)^{-1} \}$
\item\label{Odd-Facto}  $($Odd cases $m=2n+1$$)$ $L(s,\pi,\wedge^2)=\underset{k,\; i}{l.c.m.}\{ L_{ex}(s,\pi^{(2k+1)}_{i},\wedge^2)^{-1} \}$
 \end{enumerate}
 where the least common multiple is with respect to divisibility in $\mathbb{C}[q^{\pm s}]$ and is taken over all $k$ with  $k=0,1,\dotsm,n-1$ 
 and for each $k$ all constituents $\pi^{(2k)}_i$ $($ $\pi^{(2k+1)}_i$$)$ of $\pi^{(2k)}$ $($ $\pi^{(2k+1)}$, respectively$)$.
 \end{theorem}
 
 A similar definition $L_{ex}(s,\pi \times \sigma)$ and factorization formula has been constructed by Cogdell and Piatetski-Shapiro in the context of local Rankin-Selberg $L$-functions for a pair of representations $(\pi,\sigma)$ of $GL_m$ \citelist{\cite{CogPS}\cite{Matringe15}*{\S 4.1}}.

\subsection{Classifications of distinguished representations}
\label{sec:2.2}

For $m=2n$ even, we let $M_{2n}$ denote the standard Levi subgroup of $GL_{2n}$ associated with the partition $(n,n)$ of $2n$.
Let $w_{2n}=\sigma_{2n}$ and then we set $H_{2n}=w_{2n}M_{2n} w^{-1}_{2n}$. Let $w_{2n+1}=w_{2n+2}|_{GL_{2n+1}}$ so that
\[
 w_{2n+1}=\begin{pmatrix} 1 & 2 & \dotsm & n+1 & | & n+2 & n+3 & \dotsm &2n & 2n+1 \\ 
                                                1 & 3 & \dotsm & 2n+1 & | &  2 & 4 & \dotsm & 2n-2 &2n \\ 
                                                \end{pmatrix}.
\]
In the odd case, $w_{2n+1} \neq \sigma_{2n+1}$ and we denote by $M_{2n+1}$ the standard Levi subgroup attached to the partition
$(n+1,n)$ of $2n+1$. We set $H_{2n+1}=w_{2n+1}M_{2n+1}w^{-1}_{2n+1}$. We observe that $H_m$ is compatible in the sense that
$H_m \cap GL_{m-1} =H_{m-1}$. If $\alpha$ is a character of $F^{\times}$ and $\mathrm{diag}(g,g') \in M_m$, we denote by $\chi_{\alpha}$ the character 
\[
  \chi_{\alpha} : w_m \begin{pmatrix} g & \\ & g' \end{pmatrix} w^{-1}_m \mapsto \alpha \left( \frac{\mathrm{det}(g)}{\mathrm{det}(g^{\prime})} \right)
\]
of $H_m$. Let $\chi$ be a character of $H_m$. We say that $\pi \in \mathcal{A}_F(m)$ is $(H_m,\chi)$-{\it distinguished} if ${\rm Hom}_{H_m}(\pi,\chi) \neq 0$. If $\chi$ is trivial, it is customary to say that $\pi$ is  $H_m$-{\it distinguished}. In order to classify all irreducible generic distinguished representations, we need to know that the induced representations of the form ${\rm Ind}^{GL_{2n}}_Q(\Delta\otimes \widetilde{\Delta})$ are distinguished.  These types of properties over non-Archimedean local fields in characteristic zero were originally investigated by Cogdell and Piatetski-Shapiro \cite{CP94}.
Afterwards the conjecture was settled by Matringe \cite{Matringe15,Matringe17}. Parts of the proof of \cite[Proposition 3.8]{Matringe15} contain inaccuracies, and subsequently it is clarified in \cite[Proposition 5.3]{Matringe17}.

\begin{proposition}[N. Matringe] Let $\Delta$ be discrete series representations of $GL_n$ and $\alpha$ a character of $F^{\times}$. Then irreducible generic representations of the form
${\rm Ind}^{GL_{2n}}_Q(\Delta \otimes \widetilde{\Delta})$ are both $(H_{2n},\chi_{\alpha})$- and  $(S_{2n},\Theta)$-distinguished.
\end{proposition}

\begin{proof} We consider parabolically induced representations of the form $\Pi_s:={\rm Ind}^{GL_{2n}}_Q(\Delta_0\nu^s \otimes \widetilde{\Delta}_0\nu^{-s})$ with $\Delta_0$ a unitary discrete series representations of $GL_n$ and $s$ a complex parameter. The proof in \cite[Proposition 3.8]{Matringe15} relies on Bernstein's analytic continuation principle for invariant linear forms. In order to apply it to positive characteristic, we need to explain that the space ${\rm Hom}_{S_{2n}}(\Pi_s,\Theta)$ is
of dimension at most one for all $s$ except the finite number for which $\Pi_s$ is irreducible. However, if this is the case, ${\rm Hom}_{S_{2n}}(\Pi_s,\Theta)$ embeds as a subspace of ${\rm Hom}_{H_{2n} \cap P_{2n}}(\Pi_s,{\bf 1}_{H_{2n}})$ via \cite[Proposition 4.3]{Matringe} along with  ${\rm Hom}_{S_{2n}}(\Pi_s,\Theta) \subseteq {\rm Hom}_{S_{2n} \cap P_{2n}}(\Pi_s,\Theta)$. Thanks to  an auxiliary deformation parameter $s$,  the proof of \cite[Proposition 5.1-8]{Matringe15} asserts that except for a finite number of $s$,  the space  ${\rm Hom}_{H_{2n} \cap P_{2n}}(\Pi_s,{\bf 1}_{H_{2n}})$ is of dimension at most $1$, as desired. 
\par
 Alternatively, the quickest way is to use the equivalence between $(H_{2n},\chi_{\alpha})$-distinctions and $(S_{2n},\Theta)$-distinctions \cite[Corollary 3.6]{Yan22}, which only depends on Gan's approach of theta correspondence \cite[Theorem 3.1]{Gan}. This allows us to reduce to the case for  $\alpha=0$, where the result is immediate from Blanc and Delorme \cite{BD08}, as described in \cite[\S 5]{Matringe}. We refer interested reader to an excellent survey article \cite[Proposition 3.2.15]{Offen} for an expository construction of this open orbit contribution.
\end{proof}

We are now ready to introduce the classification of $(H_{2n},\chi_{\alpha})$-distinguished representations established by Matringe \cite[Theorem 3.1]{Matringe15}. 
The classification result holds 
in positive characteristic $p \neq 2$ though written in characteristic $0$ only. Indeed, the proof relies crucially on Bernstein-Zelevinsky’s version of Mackey’s theorem \cite[Theorem 5.2]{BeZe}, the explicit description of discrete series representations and their Jacquet modules \cite[Proposition 9.5]{Zelevinsky}, and the fact that a discrete series representation of $GL_{2n+1}$ cannot be $H_{2n+1}$-distinguished \cite[Theorem 3.1]{Matringe}. All the aforementioned properties are true in positive characteristic (Cf. \citelist{\cite{AKMSS}*{Appendix A}\cite{Gan}*{\S 4}}).

\begin{theorem}[N. Matringe, Even cases $m=2n$]
\label{M-Linear}
Let $\pi=\mathrm{Ind}^{GL_{2n}}_{\rm Q}(\Delta_1 \otimes \Delta_2 \otimes \dotsm \otimes \Delta_t)$ be an irreducible generic representation of $GL_{2n}$.
Let $\alpha$ be a character of $F^{\times}$ with $0 \leq \mathrm{Re}(\alpha) \leq 1/2$. Then $\pi$ is $(H_{2n},\chi_{\alpha})$-distinguished if and only if there is a reordering of the $\Delta_i$'s and an integer $r$ between $1$ and $[t/2]$,
such that $\Delta_{i+1}=\widetilde{\Delta}_i$ for $i=1,3,\dotsm,2r-1$, and $\Delta_i$ is $(H_{2n_i},\chi_{\alpha})$-distinguished for $i > 2r$.
\end{theorem}

For a discrete series representation $\Delta$ of $GL_{2n}$, $\Delta$ is $H_{2n}$-distinguished if and only if $(S_{2n},\Theta)$-distinguished. Matringe \cite[\S 5]{Matringe}, using an analytic approach, and Gan \cite[Theorem 4.2]{Gan}, using the theta correspondence, individually settled this connection. Combining this with \cite[Theorem 1.1, Proposition 5.3]{Matringe17}, we classify  $(S_{2n},\Theta)$-distinguished generic representation of $GL_{2n}$ in terms of $(S_{2n_i},\Theta)$-distinguished discrete series representations $\Delta_i$ \cite[Corollary 1.1]{Matringe17}. We refer the reader to \cite{Matringe17} for further details of the proof.

\begin{theorem}[N. Matringe]
\label{M-Shalika}
Let $\pi=\mathrm{Ind}^{GL_{2n}}_{\rm Q}(\Delta_1 \otimes \Delta_2 \otimes \dotsm \otimes \Delta_t)$ be an irreducible generic representation of $GL_{2n}$.
Then $\pi$ is $(S_{2n},\Theta)$-distinguished if and only if there is a reordering of the $\Delta_i$'s and an integer $r$ between $1$ and $[t/2]$,
such that $\Delta_{i+1}=\widetilde{\Delta}_i$ for $i=1,3,\dotsm,2r-1$, and $\Delta_i$ is $(S_{2n_i},\Theta)$-distinguished for $i > 2r$.
\end{theorem}

In the light of Theorem \ref{M-Linear}, Theorem \ref{M-Shalika}, and \cite[Theorem 4.2]{Gan}, Matringe and Gan's equivalence is valid in more general setting of irreducible generic representations of $GL_{2n}$.

\subsection{Deformations and specializations}
\label{deform-special}

Let $\pi=\mathrm{Ind}^{GL_m}_{\rm Q}(\Delta_1 \otimes \Delta_2 \otimes \dotsm \otimes \Delta_t)$ be a parabolically induced representation of $GL_m$. 
Let $\mathcal{D}_{\pi}$ denote the complex manifold $( \mathbb{C} \slash \frac{2 \pi i}{\log(q) } \mathbb{Z})^t$. The isomorphism $\mathcal{D}_{\pi} \rightarrow (\mathbb{C}^{\times})^t$
is defined by $u:=(u_1,u_2,\dotsm,u_t) \mapsto q^u:=(q^{u_1},q^{u_2},\dotsm,q^{u_t})$. We use $q^{\pm u}$ as short for $(q^{\pm u_1},q^{\pm u_2},\dotsm,q^{\pm u_t})$.
For $u \in \mathcal{D}_{\pi}$, we set
\[
 \pi_u=\mathrm{Ind}^{GL_m}_{\rm Q}(\Delta_1\nu^{u_1} \otimes \Delta_2\nu^{u_2} \otimes \dotsm \otimes \Delta_t\nu^{u_t}).
\]
Let us set
\[
  \pi^{(k_1,k_2,\dotsm,k_t)}_u=\mathrm{Ind}^{GL_m}_{\rm Q}(\Delta^{(k_1)}_1\nu^{u_1} \otimes \Delta^{(k_2)}_2\nu^{u_2} \otimes \dotsm \otimes \Delta^{(k_t)}_t\nu^{u_t}).
\]
Section \ref{deform-special} is indebted to Cogdell and Piatetski-Shapiro \cite{CogPS} and we closely follow the path of the adaptation of \cite{CogPS} used in \cite{Matringe09,Matringe15,Jo20-2} to study the characteristic zero case. In particular the deformation and specialization argument is widely available in the literature \cite{CogPS,Matringe09,Matringe15,Jo20-2}. Henceforth we only remark on the nature of the difference but the reader should consult to \cite{CogPS,Matringe09} for a complete regard.

\begin{definition}
We say that $u=(u_1,u_2,\dotsm,u_t) \in \mathcal{D}_{\pi}$ is in {\it general position} if it satisfies the following properties:
 \begin{enumerate}[label=$(\mathrm{\roman*})$]
 \item\label{General-1} For every sequences of nonnegative integers $(k_1,k_2,\dotsm,k_t)$,
 a nonzero representation
 \[
  \pi_u^{(k_1,k_2,\dotsm,k_t)}=\mathrm{Ind}(\Delta^{(k_1)}_1\nu^{u_1} \otimes \Delta^{(k_2)}_2\nu^{u_2} \otimes \dotsm \otimes  \Delta^{(k_t)}_t\nu^{u_t})
 \]
 is irreducible;
 \item\label{General-2}
 If $(a_1r_1,a_2r_2,\dotsm,a_tr_t)$ and $(b_1r_1,b_2r_2,\dotsm,b_tr_t)$ are two different sequences such that 
 \[
 \sum_{i=1}^t a_ir_i=\sum_{i=1}^t b_ir_i
 \]
 then two representations
 \[
  \mathrm{Ind}(\Delta^{(a_1r_1)}_1\nu^{u_1} \otimes \Delta^{(a_2r_2)}_2\nu^{u_2} \otimes \dotsm \otimes  \Delta^{(a_tr_t)}_t\nu^{u_t}) \quad \text{and} \quad
  \mathrm{Ind}(\Delta^{(b_1r_1)}_1\nu^{u_1} \otimes \Delta^{(b_2r_2)}_2\nu^{u_2} \otimes \dotsm \otimes \Delta^{(b_tr_t)}_t\nu^{u_t})
 \]
 possess distinct central characters;
 \item\label{General-3} If $(i,j,k,\ell) \in \{1,2,\dotsm,t \}$ with $\{i,j \} \neq \{k,\ell \}$, $L(s,\Delta_i\nu^{u_i} \times \Delta_j\nu^{u_j})$
 and $L(s,\Delta_k\nu^{u_k} \times \Delta_{\ell}\nu^{u_{\ell}})$ have no common poles;
\item\label{General-4} If $(i,j) \in \{1,2,\dotsm,t \}$ with $i \neq j$, then $L(s,\Delta_i\nu^{u_i},\wedge^2)$ and $L(s,\Delta_j\nu^{u_j},\wedge^2)$
have no common poles;
\item\label{General-5} If $(i,j,k) \in \{1,2,\dotsm,t \}$ with $i \neq j$, then $L(s,\Delta_i\nu^{u_i} \times \Delta_j\nu^{u_j})$ and $L(s,\Delta_k\nu^{u_k},\wedge^2)$
have no common poles;
\item\label{General-6} If $1 \leq i \neq j \leq t$ and $(\Delta_i^{(a_ir_i)})^{\sim} \simeq \Delta_j^{(a_jr_j)}\nu^{e}$ for some complex number $e$,
then the dimension of the space
\[
{\rm Hom}_{P_{2(n_i-a_ir_i)} \cap S_{2(n_i-a_ir_i)}}(\mathrm{Ind}(\Delta_i^{(a_ir_i)}\nu^{\frac{u_i+u_j+e}{2}}\otimes (\Delta_i^{(a_ir_i)}\nu^{\frac{u_i+u_j+e}{2}})^{\sim} ),\Theta)
\]
is at most $1$.
 \end{enumerate}
\end{definition}

We confirm that off a finite number of hyperplanes in $u$, the deformed representation $\pi_u$ is in general position \cite[Proposition 4.1]{Jo20-2}.
The important point is that $u \in \mathcal{D}_{\pi}$ in general position depends only on the representation $\pi$ not $s \in \mathbb{C}$. 
The purpose of \ref{General-2} is that outside a finite number of hyperplanes, the central character of $\pi^{(a_1r_1,a_2r_2,\dotsm,a_tr_t)}_u$ are distinct
and therefore there are only trivial extension among these representation. As a result, off these hyperplanes, the derivatives $\pi_u^{(k)}=\oplus \pi_u^{(a_1r_1,a_2r_2,\dotsm,a_tr_t)}$
are completely reducible, where $k=\sum_{i=1}^t a_ir_i$ and each $\pi_u^{(a_1r_1,a_2r_2,\dotsm,a_tr_t)}$ is irreducible. Conditions \ref{General-1} and \ref{General-2} ensure that 
 Theorem \ref{factor} is applicable. The purpose of Condition \ref{General-6} is that the occurrence of the exceptional pole of $L(s,\pi,\wedge^2)$ at $s=0$ can be determined by the existence of Shalika functional from \cite[Lemma 3.2]{Jo20-2}. Throughout Section \ref{deform-special}, we assume the Working Hypothesis proposed by E. Kaplan \cite[Remark 4.18]{Kaplan} for fields of odd characteristic.

\begin{whypothesis} Let $\Delta$ be a $(S_{2n},\Theta)$-distinguished discrete series representation of $GL_{2n}$. Then $\Delta$ is self-dual. Namely, $\widetilde{\Delta} \simeq \Delta$.
\end{whypothesis}

The following statement is a consequence of Working Hypothesis along with Theorem \ref{M-Shalika}.

\begin{corollary}
We assume the Working Hypothesis. Let $\pi=\mathrm{Ind}^{GL_{2n}}_{\rm Q}(\Delta_1 \otimes \Delta_2 \otimes \dotsm \otimes \Delta_t)$ be an irreducible generic representation of $GL_{2n}$. If $\pi$ is $(S_{2n},\Theta)$-distinguished, then $\pi$ is self-dual. Namely, $\widetilde{\pi} \simeq \pi$.
\end{corollary}

Working Hypothesis needs not to be considered for the subclass of irreducible principal series representations induced from Borel subgroups because of Theorem \ref{M-Shalika},
and we shall verify the presumption case by case in \S \ref{sec2}.

 \begin{proposition}
 \label{principal-selfdual}
 Let $\pi=\mathrm{Ind}^{GL_{2n}}_{B_{2n}}(\chi_1 \otimes \chi_2 \otimes \dotsm \otimes \chi_{2n})$ be a $(S_{2n},\Theta)$-distinguished irreducible principal series representation of $GL_{2n}$.
 Then $\pi$ is self-dual. Namely, $\widetilde{\pi} \simeq \pi$.
 \end{proposition}

Now we provide an interpretation of Theorem \ref{M-Shalika} in terms of local $L$-functions, which is analogous to Proposition 4.13 in \cite{Matringe15}.

\begin{proposition}
\label{three}
 Let $\pi=\mathrm{Ind}^{GL_{2n}}_{\rm Q}(\Delta_1 \otimes \Delta_2 \otimes \dotsm \otimes \Delta_t)$ be an irreducible generic representation of $GL_{2n}$ 
 where each $\Delta_i$ is a discrete series representation of $GL_{n_i}$ with $2n=\sum_{i=1}^t n_i$ and $t \geq 2$.
 Suppose that $L_{ex}(s,\pi,\wedge^2)$ has a pole at $s=s_0$. Then we are in the one of the following:
 \begin{enumerate}[label=$(\mathrm{\roman*})$]
\item\label{Decomp-1} There are $(i,j) \in \{1,2,\dotsm,t \}$ with $i \neq j$ such that $n_i$ and $n_j$ are even, and
 $L_{ex}(s,\Delta_i,\wedge^2)$ and  $L_{ex}(s,\Delta_j,\wedge^2)$ have $s=s_0$ as a common pole;
\item\label{Decomp-2} There are $(i,j,k,\ell) \in \{1,2,\dotsm,t \}$ with $\{i,j \} \neq \{ k,\ell \}$ such that
$L_{ex}(s,\Delta_i \times \Delta_j)$ and $L_{ex}(s,\Delta_k \times \Delta_{\ell})$ have $s=s_0$ as a common pole;
\item\label{Decomp-3} There are $(i,j,k) \in \{1,2,\dotsm,t \}$ with $i \neq j$ such that $n_k$ is even
and $L_{ex}(s,\Delta_i \times \Delta_j)$ and $L_{ex}(s,\Delta_k,\wedge^2)$ have $s=s_0$
as a common pole.
\end{enumerate}
\end{proposition}

\begin{proof}
Suppose that $L_{ex}(s,\pi,\wedge^2)$ has a pole at $s=s_0$. As $L(s,\pi,\wedge^2)=L(s-s_0,\pi\nu^{\frac{s_0}{2}},\wedge^2)$, the representation $\pi\nu^{\frac{s_0}{2}}$ admits
a non-trivial Shalika functional. We know from Theorem \ref{M-Shalika} that  $\pi\nu^{\frac{s_0}{2}}$ is isomorphic to
\[
 \mathrm{Ind}((\Delta_{i_1}\nu^{\frac{s_0}{2}}\otimes (\Delta_{i_1}\nu^{\frac{s_0}{2}})^{\sim})\otimes \dotsm \otimes (\Delta_{i_r}\nu^{\frac{s_0}{2}}\otimes (\Delta_{i_r}\nu^{\frac{s_0}{2}})^{\sim}) \otimes 
\Delta_{i_{r+1}}\nu^{\frac{s_0}{2}} \otimes \dotsm \otimes \Delta_{{i_t}}\nu^{\frac{s_0}{2}} ),\; 0 \leq r \leq [t/2],
\]
where $\Delta_{i_j}\nu^{\frac{s_0}{2}}$ affords a Shalika functional and each $n_{i_j}$ is even for all $j > r$. Putting it in a different way, $(\Delta_i\nu^{\frac{s_0}{2}})^{\sim} \simeq \Delta_j\nu^{\frac{s_0}{2}}$ with $i \neq j$, or $\Delta_k\nu^{\frac{s_0}{2}}$ owns a Shalika functional with $n_k$ an even number.

\par
According to \cite[Proposition 4.6]{Matringe15}, $(\Delta_i\nu^{\frac{s_0}{2}})^{\sim} \simeq \Delta_j\nu^{\frac{s_0}{2}}$ or equivalently $\widetilde{\Delta}_i \simeq \Delta_j\nu^{s_0}$
if and only if $L_{ex}(s,\Delta_i \times \Delta_j)$ has a pole at $s=s_0$.

\par
If $\Delta_k\nu^{\frac{s_0}{2}}$ has the Shalika functional, the space ${\rm Hom}_{S_{n_k}}(\Delta_k\nu^{\frac{s_0}{2}},\Theta)$ is non-trivial and that its central character $\omega_{\Delta_k\nu^{s_0/2}}$ is trivial. Since $\Delta_k\nu^{\frac{s_0}{2}}$ is the irreducible square integrable representation, we obtain from \cite[Proposition 3.4]{Jo20-2}
that $L_{ex}(s,\Delta_k\nu^{\frac{s_0}{2}},\wedge^2)$ has a pole at $s=0$ or equivalently $L_{ex}(s,\Delta_k,\wedge^2)$ has a pole at $s=s_0$. Therefore $s=s_0$ is the 
common pole for either of three cases in Proposition \ref{three}.
\end{proof}

Let $\Delta$ be a discrete series representation. According to \cite{Zelevinsky}, such a representation $\Delta$ is the unique irreducible quotient of the form:
$\mathrm{Ind}^{GL_m}_{\rm Q}(\rho \otimes \rho\nu \otimes \dotsm \otimes \rho\nu^{\ell-1})$, 
where the induction is normalized parabolic induction from the standard parabolic subgroup $\rm Q$ attached to the partition $(r,r,\dotsm,r)$ of $m=r\ell$
and $\rho \in \mathcal{A}_F(r)$ is irreducible and supercuspidal. We denote by $\Delta=[\rho,\rho\nu,\dotsm,\rho\nu^{\ell-1}]$ such a quotient.
Using Hartogs' Theorem \cite{Jo20-2} is more close to the original spirit of the direction in \cite{CogPS}. Nevertheless we present an alternative approach employing Proposition \ref{three} to keep the uniformity with \cite{Matringe09,Matringe15}.

\begin{proposition}
\label{deform-L}
We assume the Working Hypothesis. Let $\pi=\mathrm{Ind}^{GL_m}_{\rm Q}(\Delta_1 \otimes \Delta_2 \otimes \dotsm \otimes \Delta_t)$ be a parabolically induced representation of $GL_m$.
Let $u=(u_1,u_2,\dotsm,u_t) \in \mathcal{D}_{\pi}$ be in general position and $\pi_u=\mathrm{Ind}^{GL_m}_{\rm Q}(\Delta_1\nu^{u_1} \otimes \Delta_2\nu^{u_2} \otimes \dotsm \otimes \Delta_t\nu^{u_t})$
the deformed representation. Then we have the following:
 \begin{enumerate}[label=$(\mathrm{\roman*})$]
\item\label{Deform-L-1} $\displaystyle L(s,\pi_u,\wedge^2)=\prod_{1 \leq k \leq t} L(s+2u_k,\Delta_k,\wedge^2) \prod_{1 \leq i < j  \leq t} L(s+u_i+u_j, \Delta_i \times \Delta_j)$.
\item\label{Deform-L-2} There is a polynomial $Q(X) \in \mathbb{C}[X]$ such that
\[
L(s,\pi,\wedge^2)
=Q(q^{-s})\prod_{1 \leq k \leq t} L(s,\Delta_k,\wedge^2) \prod_{1 \leq i < j  \leq t} L(s, \Delta_i \times \Delta_j).
\]
 \end{enumerate}
\end{proposition}

\begin{proof}
Let us take $\Delta_i$ to be associated to the segment $[\rho_i,\rho_i\nu,\dotsm,\rho_i\nu^{\ell_i-1}]$ with $\rho_i$ an irreducible supercuspidal representation of $GL_{r_i}$, $m_i=r_i\ell_i$, and $m=\sum_{i=1}^tr_i\ell_i$. Keeping Theorem \ref{factor} in hand, we set out
\[
 L(s,\pi_u,\wedge^2)^{-1}=l.c.m. \{ L_{ex}(s,\pi^{(a_1r_1,a_2r_2,\dotsm,a_tr_t)}_u,\wedge^2)^{-1} \}
\]
where $0 \leq a_i \leq \ell_i$, $m-\sum_{i=1}^t a_ir_i$ is an even number, and the least common multiple is taken  
in terms of divisibility in $\mathbb{C}[q^{\pm s}]$. Suppose that $L_{ex}(s,\pi^{(a_1r_1,a_2r_2,\dotsm,a_tr_t)}_u,\wedge^2)$ has a pole at $s=s_0$. If the number of indices $i$
such that $r_i \neq \ell_i$ is more than $3$, we deduce from Proposition \ref{three} that
 \begin{enumerate}[label=$(\mathrm{\roman*})$]
\item\label{Excep-1} There are $(i,j) \in \{1,2, \dotsm,t\}$ with $i \neq j$ such that $m_i-a_ir_i$ and $m_j-a_jr_j$ are even, and
$L(s,\Delta_i^{(a_ir_i)}\nu^{u_i},\wedge^2)$ and $L(s,\Delta_j^{(a_jr_j)}\nu^{u_j},\wedge^2)$ have $s=s_0$ as a common pole;
\item\label{Excep-2} There are $(i,j,k,\ell) \in \{1,2, \dotsm,t\}$ with $\{i,j \} \neq \{k,\ell \}$ such that $L(s,\Delta_i^{(a_ir_i)}\nu^{u_i} \times \Delta_j^{(a_jr_j)}\nu^{u_j})$
and $L(s,\Delta_k^{(a_kr_k)}\nu^{u_k} \times \Delta_{\ell}^{(a_{\ell}r_{\ell})}\nu^{u_{\ell}})$ have $s=s_0$ as a common pole;
\item\label{Excep-3} There are $(i,j,k) \in \{1,2, \dotsm,t\}$ with $i \neq j$ such that $m_k-a_kr_k$ is even, and
$L(s,\Delta_i^{(a_ir_i)}\nu^{u_i} \times \Delta_j^{(a_jr_j)}\nu^{u_j})$ and $L(s,\Delta_k^{(a_kr_k)}\nu^{u_k},\wedge^2)$
have $s=s_0$ as a common pole.
\end{enumerate}
However, Conditions \ref{General-3}, \ref{General-4}, and \ref{General-5} of general positions ensure that the above scenario
cannot happen as long as $u$ is in general position, because exceptional poles are poles of original $L$-factors
$L(s,\Delta_i \times \Delta_j)$ and $L(s,\Delta_k,\wedge^2)$. Owing to \cite[Corollary 4.11]{Jo20-2}, when there exists exactly one pair $(i,j)$
of indices $i \neq j$ such that $r_i \neq \ell_i$ and $r_j \neq \ell_j$, we have
\[
L_{ex}(s,\mathrm{Ind}(\Delta_i^{(a_ir_i)}\nu^{u_i} \otimes \Delta_j^{(a_jr_j)}\nu^{u_j}),\wedge^2)=L_{ex}(s,\Delta_i^{(a_ir_i)}\nu^{u_i}  \times  \Delta_j^{(a_jr_j)}\nu^{u_j}).
\]
if $i$ is the only index such that $r_i \neq \ell_i$, it is nothing but $L_{ex}(s,\Delta_i^{(a_ir_i)}\nu^{u_i},\wedge^2)$.

\par
To summarize $L_{ex}(s,\pi^{(a_1r_1,a_2r_2,\dotsm,a_tr_t)}_u,\wedge^2)$ is equal to either $L_{ex}(s,\Delta_i^{(a_ir_i)}\nu^{u_i}  \times  \Delta_j^{(a_jr_j)}\nu^{u_j})$ for $i < j$  or $L_{ex}(s,\Delta_i^{(a_ir_i)}\nu^{u_i},\wedge^2)$. Following the rest of the proof in \cite[Theorem 5.1]{Jo20-2}, we arrive at
\[
  \begin{split}
  L(s,\pi_u,\wedge^2)&=\prod_{1 \leq k \leq t} L(s,\Delta_k\nu^{u_k},\wedge^2) \prod_{1 \leq i < j \leq t} L(s,\Delta_i\nu^{u_i}\times \Delta_j\nu^{u_j} )\\
  &=\prod_{1 \leq k \leq t} L(s+2u_k,\Delta_k,\wedge^2) \prod_{1 \leq i < j \leq t} L(s+u_i+u_j,\Delta_i\times \Delta_j ).
  \end{split}
\]

\par
Concerning the second part, let $\mathcal{W}_{\pi}^{(0)}$ be the Whittaker model associated to $\pi_u$ \cite[\S 3.1]{CogPS}. For $W_u \in \mathcal{W}_{\pi}^{(0)}$, it follows from 
the standard Bernstein's principle of meromorphic continuation and rationality \cite[Proposition 4.2, Proposition 4.4]{Jo20-2} that $J(s,W_u,\Phi)$ defines a rational function in $\mathbb{C}(q^{-s},q^{-u})$. We conclude the first part \ref{Deform-L-1} that the rational function
\[
  \frac{J(s,W_u,\Phi)}{\prod_{1 \leq k \leq t} L(s+2u_k,\Delta_k,\wedge^2) \prod_{1 \leq i < j \leq t} L(s+u_i+u_j,\Delta_i\times \Delta_j )}
\]
has no poles on the Zariski open set of $u$ in general position. We can take one step further to assert that the ratio lies in $\mathbb{C}[q^{\pm s},q^{\pm u}]$
by the proof of \cite[Lemma 5.1]{Matringe15} and \cite[Proposition 5.3]{Jo20-2}. The statement is now an immediate consequence of
specialization to $u=0$.
\end{proof}

We denote by $P \sim Q$ that the ratio is a unit in $\mathbb{C}[q^{\pm s}]$ for two rational functions $P(q^{-s})$ and $Q(q^{-s})$ in $\mathbb{C}(q^{-s})$. As alluded in Langlands-Shahidi method \cite{GL,HL-Exterior,HL-RS,Lomeli}, the unit emerging in Theorem \ref{Deform-gamma}-\ref{Deform-gamma-2} will be presumably $1$. This is so-called the 
{\it multiplicativity of} $\gamma$-{\it factors}. However demonstrating the multiplicativity property requires manipulating integrals in a delicate manner. Nonetheless, it seems likely that the weaker one that is relevant to us is enough for the application therein.

\begin{theorem}
\label{Deform-gamma}
We assume the Working Hypothesis. Let $\pi=\mathrm{Ind}^{GL_m}_{\rm Q}(\Delta_1 \otimes \Delta_2 \otimes \dotsm \otimes \Delta_t)$ be a parabolically induced representation of $GL_m$.
Let $u=(u_1,u_2,\dotsm,u_t) \in \mathcal{D}_{\pi}$ be in general position and $\pi_u=\mathrm{Ind}^{GL_m}_{\rm Q}(\Delta_1\nu^{u_1} \otimes \Delta_2\nu^{u_2} \otimes \dotsm \otimes \Delta_t\nu^{u_t})$
the deformed representation. Then we have the following:
 \begin{enumerate}[label=$(\mathrm{\roman*})$]
\item\label{Deform-gamma-1} $\displaystyle \gamma(s,\pi_u,\wedge^2,\psi) \sim \prod_{1 \leq k \leq t} \gamma(s+2u_k,\Delta_k,\wedge^2,\psi) \prod_{1 \leq i < j  \leq t} \gamma(s+u_i+u_j, \Delta_i \times \Delta_j,\psi)$.
\item\label{Deform-gamma-2} $\displaystyle \gamma(s,\pi,\wedge^2,\psi)
\sim \prod_{1 \leq k \leq t} \gamma(s,\Delta_k,\wedge^2,\psi) \prod_{1 \leq i < j  \leq t} \gamma(s, \Delta_i \times \Delta_j,\psi)$.
\end{enumerate} 
\end{theorem}

\begin{proof}
The proof proceeds along the line of \cite[Proposition 5.4]{Jo20-2} and \cite[Proposition 5.5]{Matringe15} by appealing Theorem \ref{exterior-func} and  Proposition \ref{deform-L} to our framework and it is originated from Cogdell and Piatetski-Shapiro \cite[Proposition 4.3]{CogPS}. The second statement \ref{Deform-gamma-2} can be acomplished by specializing to $u=0$.
\end{proof}

To proceed further, we adopt the terminology from \cite{CogPS,Matringe15}. We say that $\pi \in \mathcal{A}_F(m)$ is a representation of 
{\it Langlands type} if $\Xi$ is of the form $\mathrm{Ind}^{GL_m}_{\rm Q}({\Delta_{\circ}}_1\nu^{u_1} \otimes {\Delta_{\circ}}_2\nu^{u_2} \otimes \dotsm \otimes {\Delta_{\circ}}_t\nu^{u_t})$,
where each ${\Delta_{\circ}}_i$ is the irreducible square integrable representation of $GL_{m_i}$, $m_1+m_2+\dotsm+m_t=m$, each $u_i$ is real and they are 
ordered so that $u_1 \geq u_2 \geq \dotsm \geq u_t$. Let $\pi$ be an irreducible admissible representation of $GL_m$. Regardless of being generic, 
$\pi$ can be realized as the unique Langlands quotient of Langlands type $\Xi=\mathrm{Ind}^{GL_m}_{\rm Q}({\Delta_{\circ}}_1\nu^{u_1} \otimes {\Delta_{\circ}}_2\nu^{u_2} \otimes \dotsm \otimes {\Delta_{\circ}}_t\nu^{u_t})$ which is of Whittaker type. The exterior square $L$-factor is defined to be
\[
  L(s,\Xi,\wedge^2)=L(s,\pi,\wedge^2).
\]

\begin{theorem}
\label{ext-langlands-prod}
We assume the Working Hypothesis. Let $\pi=\mathrm{Ind}^{GL_m}_{\rm Q}({\Delta_{\circ}}_1\nu^{u_1} \otimes {\Delta_{\circ}}_2\nu^{u_2} \otimes \dotsm \otimes {\Delta_{\circ}}_t\nu^{u_t})$ be
a representation of Langlands type of $GL_m$. Then we have
\[
  L(s,\pi,\wedge^2)=\prod_{1 \leq k \leq t} L(s+2u_k,{\Delta_{\circ}}_k,\wedge^2) \prod_{1 \leq i < j  \leq t} L(s+u_i+u_j, {\Delta_{\circ}}_i \times {\Delta_{\circ}}_j).
\]
\end{theorem}

\begin{proof}
The proof is akin to those of \cite[Theorem 4.1]{CogPS}, \cite[Theorem 5.7]{Jo20-2} and \cite[Theorem 4.26]{Matringe09}. In order to be concise, we do not include the complete details.
\end{proof}

We pass to the case of irreducible generic representations.

\begin{corollary}
\label{generic}
We assume the Working Hypothesis. Let $\pi=\mathrm{Ind}^{GL_m}_{\rm Q}(\Delta_1 \otimes \Delta_2 \otimes \dotsm \otimes \Delta_t)$ be an irreducible generic representation of $GL_m$. Then we have
\[
  L(s,\pi,\wedge^2)=\prod_{1 \leq k \leq t} L(s,\Delta_k,\wedge^2) \prod_{1 \leq i < j  \leq t} L(s, \Delta_i \times \Delta_j).
\]
\end{corollary}

\begin{proof}
Since $\pi$ is irreducible, essentially square integrable representations $\Delta_i$ can be rearranged to be in Langlands order without changing $\pi$.
\end{proof}

We define the symmetric square $L$-factor to be the ratio of Rankin-Selberg $L$-factors for $GL_m \times GL_m$ by exterior square $L$-factors for $GL_m$:
\begin{equation}
\label{ext-sym}
 L(s,\pi,\mathrm{Sym}^2)=\frac{L(s,\pi \times \pi)}{  L(s,\pi,\wedge^2)}.
\end{equation}
In comparison to \cite{Matringe09} and \cite{Matringe15}, we pursue purely local means more to express 
a local exterior square $L$-function in terms of local $L$-functions for supercuspidal representations.
Performing this step has benefits to make the globalization result of Henniart and Lomel\'{\i} \cite{HL-Exterior,HL-RS} feasible, instead of globalizing discrete series representations \cite{Kaplan,KeRa,Matringe09} as a black box.

\begin{theorem}
\label{Discrete}
We assume that the Working Hypothesis holds for the subclass of all irreducible supercuspidal representations. Let $\Delta_{\circ}=[\rho_{\circ}\nu^{-\frac{\ell-1}{2}},\dotsm,\rho_{\circ}\nu^{\frac{\ell-1}{2}}]$ be an irreducible square integrable representation of $GL_{\ell r}$
with $\rho_{\circ}$ an irreducible unitary supercuspidal representation of $GL_r$.
 \begin{enumerate}[label=$(\mathrm{\roman*})$]
\item\label{Discrete-1} Suppose that $\ell$ is even. Then we have
\[
\begin{split}
 L(s,\Delta_{\circ},\wedge^2)=&\prod_{i=1}^{\ell/2}L(s,\rho_{\circ}\nu^{(\ell+1)/2-i},\wedge^2)L(s,\rho_{\circ}\nu^{\ell/2-i},\mathrm{Sym}^2);\\
 L(s,\Delta_{\circ},\mathrm{Sym}^2)=&\prod_{i=1}^{\ell/2}L(s,\rho_{\circ}\nu^{(\ell+1)/2-i},\mathrm{Sym}^2)L(s,\rho_{\circ}\nu^{\ell/2-i},\wedge^2).\\
\end{split}
\]
\item\label{Discrete-2} Suppose that $\ell$ is odd. Then we have
\[
\begin{split}
 L(s,\Delta_{\circ},\wedge^2)=&\prod_{i=1}^{(\ell+1)/2}L(s,\rho_{\circ}\nu^{(\ell+1)/2-i},\wedge^2)\prod_{i=1}^{(\ell-1)/2}L(s,\rho_{\circ}\nu^{\ell/2-i},\mathrm{Sym}^2);\\
 L(s,\Delta_{\circ},\mathrm{Sym}^2)=&\prod_{i=1}^{(\ell+1)/2}L(s,\rho_{\circ}\nu^{(\ell+1)/2-i},\mathrm{Sym}^2)\prod_{i=1}^{(\ell-1)/2}L(s,\rho_{\circ}\nu^{\ell/2-i},\wedge^2);\\
\end{split}
\]
\end{enumerate}
\end{theorem}

\begin{proof} Our proof is truly influenced by Shahidi \cite[Proposition 8.1]{Shahidi}. 
By the uniqueness of the Whittaker functional, the Whittaker model for $\Delta_{\circ}$ agrees with that for $\xi=\mathrm{Ind}(\rho_{\circ}\nu^{-\frac{\ell-1}{2}} \otimes \dotsm \otimes \rho_{\circ}\nu^{\frac{\ell-1}{2}})$. Likewise the same feature holds for $\xi^{\iota}:=\mathrm{Ind}(\widetilde{\rho}_{\circ}\nu^{-\frac{\ell-1}{2}}\otimes \dotsm \otimes \widetilde{\rho}_{\circ}\nu^{\frac{\ell-1}{2}})$ and $\widetilde{\Delta}_{\circ}$. This puts us in a position to manifest that
\[
 \gamma(s,\Delta_{\circ},\wedge^2,\psi)=\gamma(s,\mathrm{Ind}(\rho_{\circ}\nu^{-\frac{\ell-1}{2}} \otimes \dotsm \otimes \rho_{\circ}\nu^{\frac{\ell-1}{2}}),\wedge^2,\psi).
\]
Let $u$ be in general position and $\xi_u=\mathrm{Ind}(\rho_{\circ}\nu^{u_1-\frac{\ell-1}{2}} \otimes  \dotsm \otimes \rho_{\circ}\nu^{u_{\ell}+\frac{\ell-1}{2}})$ its associated deformed representation. Upon noting the assumption that any $(S_{2n},\Theta)$-distinguished irreducible supercuspidal representation $\rho$ is self-dual, we see that Proposition \ref{three} to Theorem \ref{Deform-gamma} can be completely carried over verbatim to the triple $(\Delta_{\circ},\xi,\xi_u)$. The remainder of the proof is parallel to that of \cite[Theorem 5.12]{Jo20-2} (cf. proof of Proposition \ref{Discrete-BF}),
 and we find
\[
 L(s,\Delta_{\circ},\wedge^2)=
  \begin{cases}
 \displaystyle \prod_{i=1}^{\ell/2}L(s,\rho_{\circ}\nu^{(\ell+1)/2-i},\wedge^2)L(s,\rho_{\circ}\nu^{\ell/2-i},\mathrm{Sym}^2) & \text{if $\ell$ is even,} \\
\displaystyle  \prod_{i=1}^{(\ell+1)/2}L(s,\rho_{\circ}\nu^{(\ell+1)/2-i},\wedge^2)\prod_{i=1}^{(\ell-1)/2}L(s,\rho_{\circ}\nu^{\ell/2-i},\mathrm{Sym}^2) &  \text{if $\ell$ is odd.} \\
 \end{cases}
\]
The expression of local symmetric square $L$-function $L(s,\Delta_{\circ},{\rm Sym}^2)$ is a direct consequence of the factorization $L(s,\Delta_{\circ} \times \Delta_{\circ})=L(s,\Delta_{\circ},\wedge^2)L(s,\Delta_{\circ},{\rm Sym}^2)$ just as in \eqref{ext-sym}.
\end{proof}

\subsection{The equality for principal series representations}

We briefly review the Langlands-Shahidi method to the local exterior square $L$-function \cite{GL,HL-Exterior}. Let $\mathbf{G}=Sp_{2m}$ be a symplectic group over $F$ in $2m$ variables. The group $\mathbf{M} \simeq GL_m$ can be embedded as a Levi component of a maximal Siegel parabolic subgroup $\mathbf{P}=\mathbf{M}\mathbf{N}$ with unipotent radical $\mathbf{N}$. Let $r$ be the adjoint representation of the $L$-group of $\mathbf{M}$ on ${^L\mathfrak{n}}$, the Lie algebra of the $L$-group of $\mathbf{N}$.
We can check that $r=r_1 \oplus r_2$. The irreducible representation $r_1$ gives the standard $\gamma$-factor of $GL_n$ and $r_2$ gives the {\it Langlands-Shahidi exterior square} $\gamma$-{\it factor}, $\gamma(s,\pi,r_2,\psi)=\gamma_{LS}(s,\pi,\wedge^2,\psi)$. The $\gamma$-factor $\gamma_{LS}(s,\pi,\wedge^2,\psi)$ defined in \cite{HL-Exterior} is a rational function in $\mathbb{C}(q^{-s})$. Let $P(X)$ be the unique polynomial in $\mathbb{C}[X]$ satisfying $P(0)=1$ and such that $P(q^{-s})$ is the numerator of $\gamma_{LS}(s,\pi,\wedge^2,\psi)$. Whenever $\pi$ is tempered, the {\it local Langlands-Shahidi exterior square} $L$-{\it function} is defined by
\[
 \mathcal{L}(s,\rho,\wedge^2):=P(q^{-s})^{-1}.
\]
We observe that $\pi$ tempered implies that $\mathcal{L}(s,\rho,\wedge^2)$ is holomorphic for $\mathrm{Re}(s) > 0$ \cite[\S 4.6]{HL-Exterior}. The {\it Langlands-Shahidi exterior square} $\varepsilon$-{\it factor} is defined to satisfy the relation:
\[
 \varepsilon_{LS}(s,\pi,\wedge^2,\psi)=\gamma_{LS}(s,\pi,\wedge^2,\psi) \frac{\mathcal{L}(1-s,\widetilde{\pi},\wedge^2)}{\mathcal{L}(s,\pi,\wedge^2)}.
\]
Besides, various types of $L$-factors $\mathcal{L}(s,\pi,\mathrm{Sym}^2)$ for $\mathbf{G}=SO_{2m+1}$, $\mathcal{L}(s,\pi \times \pi)$ for $\mathbf{G}=GL_{2m}$, and $\mathcal{L}(s,\pi,As)$ for $\mathbf{G}=U_{m}$, can be extracted from \cite{HL-RS,Lomeli}.

\begin{proposition}
\label{principal-disc}
Let $\Delta$ be a discrete series representation of the form 
\begin{equation}
\label{principal-disc-equn}
[\chi,\chi\nu,\dotsm,\chi\nu^{\ell-1}],
\end{equation}
where $\chi$ is a character of $F^{\times}$. Then we have
\[
L(s,\Delta,\wedge^2)=\mathcal{L}(s,\Delta,\wedge^2).
\]
As a consequence, if $\ell=2n$ is even and $\Delta$ is $(S_{2n},\Theta)$-distinguished, then $\Delta$ is self-dual.
\end{proposition}

\begin{proof}
Just as observed in Proposition \ref{principal-selfdual}, the Working Hypothesis does not need to be checked for the character $\chi$ of $F^{\times}$, and $\Delta$ is automatically self-dual.
As in the proof of Theorem \ref{BF-disc}, we can easily reduce it to the case with $\Delta$ a unitary representation. We are then left with applying Theorem \ref{Discrete} to $\Delta$,
 from which the equality shall follow by comparing it with the work of Shahidi \cite[Proposition 8.1]{Shahidi}.
\end{proof}

Let us turn our attention to the subclass of irreducible generic subquotients of principal series representations. This class is not necessarily spherical.

\begin{proposition}
\label{principal}
 Let $\pi$ be an irreducible generic subquotient of a principal series representation of $GL_m$. Then we have
\[
 L(s,\pi,\wedge^2)=\mathcal{L}(s,\pi,\wedge^2).
\]
\end{proposition}

\begin{proof}
From \cite{BeZe,Zelevinsky}, $\pi$ is of the form $\mathrm{Ind}(\Delta_1 \otimes \Delta_2 \otimes \dotsm \otimes \Delta_t)$, where each $\Delta_i$ is either a character $\chi_i$ of $F^{\times}$ or  a discrete series representation given by the segment of the form \eqref{principal-disc-equn}. In considering Proposition \ref{principal-disc}, any $(S_{2{n_i}},\Theta)$-distinguished representations $\Delta_i$'s satisfy the Working Hypothesis. The inductive relation formula, namely, Corollary \ref{generic} is applicable and it can be accomplished that
\[
  L(s,\pi,\wedge^2)=\prod_{1 \leq k \leq t} L(s,\Delta_k,\wedge^2) \prod_{1 \leq i < j  \leq t} L(s, \Delta_i \times \Delta_j).
\]
In the aspect of Proposition \ref{principal-disc}, we only need to compare it with \cite[Theorem 3.1.(xi)]{GL}.
\end{proof}

The unramified character $\chi$ means that it is invariant under the maximal compact subgroup $\mathcal{O}^{\times}$ of $F^{\times}$. As before, the Working Hypothesis is no longer needed for the set of irreducible unramified representations. Hence Corollary \ref{generic} in the preceding section \S \ref{deform-special}, has the following result.

\begin{corollary}
\label{unramified}
 Let $\pi=\mathrm{Ind}^{GL_m}_{B_m}(\chi_1 \otimes \chi_2 \otimes \dotsm \otimes \chi_m)$ be an irreducible full induced representation from the Borel subgroup of unramified character $\chi_i$ of $F^{\times}$. Then we have
 \[
   L(s,\pi,\wedge^2)=\prod_{1 \leq i < j \leq m}\frac{1}{1-\chi_i(\varpi)\chi_j(\varpi)q^{-s}}.
 \]
\end{corollary}

\section{Local to Global Argument}
\label{sec3}
\subsection{Eulerian integral representations}
We denote by $\mathbb{F}_q$ the residue field of $F$, and
let $k=\mathbb{F}_q(t)$ be a (global) function field of the projective line $\mathbb{P}^1$ over $\mathbb{F}_q$. 
Let $\mathbb{A}$ denote its ring of ad\`{e}les. Let $(\Pi,V_{\Pi})$ be a cuspidal automorphic representation of $GL_m(\mathbb{A})$.
We denote by $|\mathbb{P}^1|$ the set of closed points of $\mathbb{P}^1$. The set $|\mathbb{P}^1|$ is in bijection with the set of places of $k$. 
Hence we write by abuse of notation $|\mathbb{P}^1|$ for the set of places of $k$. 
Since $\Pi$ is irreducible, we have restricted tensor product decomposition $\Pi=\bigotimes'_v \Pi_v$ with $(\Pi_v,V_{\Pi_v})$ irreducible admissible generic representations of $GL_m(k_v)$ \cite{Flath} (Cf. \cite[\S 4]{Cogdell-DOC}). Let its central character be $\omega_{\Pi}$. We let $P_{n-1,1}=Z_nP_n$ be the standard parabolic subgroup associated to the partition $(n-1,1)$ of $n$. Each $\Phi \in \mathcal{S}(\mathbb{A}^n)$
defines a smooth function on $GL_n(\mathbb{A})$, left invariant by $P_n(\mathbb{A})$, by $g \mapsto \Phi(e_ng)$ for $g \in GL_n(\mathbb{A})$. We consider the function
\[
 f(s,g;\Phi,\omega_{\Pi})=|\mathrm{det}(g)|^s \int_{\mathbb{A}^{\times}} \omega_{\Pi}(z) \Phi(ze_ng)|z|^{ns} d^{\times} z
\]
with the absolute convergence of the integral \cite[(4.1)]{JaSh81}. We extend $\omega_{\pi}$ to a character of $P_{n-1,1}$ by $\omega_{\Pi}(p)=\omega_{\Pi}(a)$ for $p=\begin{pmatrix} h&u \\ & a\end{pmatrix} \in P_{n-1,1}$. We construct the Eisenstein series by
\[
 E(s,g;\Phi,\omega_{\Pi})=\sum_{\gamma \in P_{n-1,1}(k) \backslash GL_n(k)} F(s,\gamma g; \Phi,\omega_{\Pi})
\]
This series is convergent absolutely for $\mathrm{Re}(s) > 1$ \cite[(4.1)]{JaSh81}. The mirabolic (Godement-Jacquet) Eisenstein series $E(s,g;\Phi,\omega_{\Pi})$ has a meromorphic continuation to all of $\mathbb{C}$ and satisfies the following functional equation \cite[\S 4]{JaSh81}:
\begin{equation}
\label{func-eisenstein}
  E(s,g;\Phi,\omega_{\Pi})=E(1-s,{^{\iota}g};\hat{\Phi},\omega^{-1}_{\Pi})
\end{equation}
where $^{\iota}g={^tg^{-1}}$ and the Fourier transform on $\mathcal{S}(\mathbb{A}^n)$ is defined by
\[
 \hat{\Phi}(y)=\int_{\mathbb{A}^n} \Phi(x)\psi(x {\,^ty}) dx.
\]
For $m=2n$, $\Phi \in \mathcal{S}(\mathbb{A}^n)$, and $\varphi \in V_{\Pi}$, we let 
\[
\begin{split}
 &I_{\psi}(s,\varphi,\Phi)\\
 &=\int_{Z_n(k)GL_n(k) \backslash GL_n(\mathbb{A})} \int_{\mathcal{M}_n(k) \backslash \mathcal{M}_n(\mathbb{A})} \varphi \left( \begin{pmatrix} I_n & X \\ & I_n \end{pmatrix}  \begin{pmatrix} g & \\ & g \end{pmatrix}\right) \psi^{-1}(\mathrm{Tr}(X)) E(s,g:\Phi,\omega_{\Pi}) dXdg.
\end{split}
\]
For $m=2n+1$, $\Phi \in \mathcal{S}(\mathbb{A}^n)$, and $\varphi \in V_{\Pi}$, we define a global integral as follows:
\[
\begin{split}
 I_{\psi}(s,\varphi,\Phi)
 &=\int\limits_{\mathbb{A}^n}\int\limits_{GL_n(k) \backslash GL_n(\mathbb{A})} \int\limits_{\mathcal{M}_n(k) \backslash \mathcal{M}_n(\mathbb{A})} \int\limits_{k^n \backslash \mathbb{A}^n} \varphi 
 \left( \begin{pmatrix} I_n & X&Z \\ &I_n& \\ &&1  \end{pmatrix}   \begin{pmatrix} g && \\ &g& \\ & &1 \end{pmatrix}  \begin{pmatrix} I_n && \\ &I_n& \\ &y &1 \end{pmatrix} \right)\\
 &\quad \times \psi^{-1}(\mathrm{Tr}(X)) \Phi(y) |\mathrm{det}(g)|^{s-1} dZ dX dg dy.
\end{split}
\]
The following theorem gives a meaning to these global integrals.

\begin{theorem}
\label{global-func} 
The integral $I_{\psi}(s,\varphi,\Phi)$ is convergent for $\mathrm{Re}(s)$ large enough, represents a meromorphic function on the 
entire plane and satisfies the functional equation
\[
 I_{\psi}(s,\varphi,\Phi)= I_{\psi^{-1}}(1-s,\varrho(\tau_m)\widetilde{\varphi},\hat{\Phi}),
\]
where $\varrho$ denotes right translation and $\widetilde{\varphi}(g)=\varphi(^{\iota}g)$.
\end{theorem}

\begin{proof}
The analytic properties have been established for the even case $m=2n$ in Section 5 of \cite{JaSh88} and the odd case $m=2n+1$ in Section 9 of \cite{JaSh88}.
The functional equation for the even case $m=2n$ follows immediately from that of the Eisenstein series $E(s,g:\Phi,\omega_{\Pi})$ \eqref{func-eisenstein}. See also \cite[Theorem 3.11]{KeRa}.
We take this occasion to refine the elaboration for the odd case $m=2n+1$ in \cite[\S 3.5]{CoMa} thoroughly. 
If $\varphi \in V_{\Pi}$, then $\varphi_1$ and $\varphi_2$ are defined on the page 219 of \cite{JaSh88}:
\[
 \varphi_1(g)=\int_{\mathbb{A}^n} \varphi \left( g \begin{pmatrix} I_n && \\ &I_n& \\ &y &1 \end{pmatrix} \right) \Phi(y) dy; \quad
  \varphi_2(g)=\int_{\mathbb{A}^n} \varphi \left( g \begin{pmatrix} I_n &&y \\ &I_n& \\ & &1 \end{pmatrix} \right) \hat{\Phi}(-{^ty}) dy
\]
where $\Phi \in \mathcal{S}(\mathbb{A}^n)$. We begin to deal with the equation on the bottom of page 219 of \cite{JaSh88}:
\[
\begin{split}
 &\int_{k^n \backslash \mathbb{A}^n} \int_{\mathcal{M}_n(k) \backslash \mathcal{M}_n(\mathbb{A}) } \varphi_1 \left(  \begin{pmatrix} I_n &X& \\ &I_n& \\ & &1 \end{pmatrix}
 \begin{pmatrix} I_n &&Z \\ &I_n& \\ & &1 \end{pmatrix}  \begin{pmatrix} g && \\ &g& \\ & &1 \end{pmatrix} \right) \psi^{-1}(\mathrm{Tr}(X))dX dZ \\
 &=\int_{k^n \backslash \mathbb{A}^n} \int_{\mathcal{M}_n(k) \backslash \mathcal{M}_n(\mathbb{A}) } \varphi_2 \left(  \begin{pmatrix} I_n &X& \\ &I_n& \\ & &1 \end{pmatrix}
 \begin{pmatrix} I_n && \\ &I_n& \\ &Z &1 \end{pmatrix}  \begin{pmatrix} g && \\ &g& \\ & &1 \end{pmatrix} \right) \psi^{-1}(\mathrm{Tr}(X))dX dZ |\mathrm{det}(g)|.
\end{split}
\]
(Here $\varphi$ in the corresponding formula in \cite[P. 219]{JaSh88} seems to be $\varphi_2$). As opposed to Jacquet and Shalika who conjugate them with the permutation matrix $\begin{pmatrix} &w_n& \\ w_n&& \\ &&1  \end{pmatrix}$, we exploit $\tau_{2n+1}$. This articulation is consistent with the shape of the local functional equation in \cite[Theorem 3.1]{CoMa}.
By applying $g \mapsto \tau_{2n+1}{^{\iota}g}\tau^{-1}_{2n+1}$, and then changing the variables $X \mapsto -X$ and $Z \mapsto -Z$, the above integral is written as
\[
\begin{split}
 \int_{k^n \backslash \mathbb{A}^n} \int_{\mathcal{M}_n(k) \backslash \mathcal{M}_n(\mathbb{A}) } \widetilde{\varphi}_2 \left(  \begin{pmatrix} I_n &X& \\ &I_n& \\ & &1 \end{pmatrix}
 \begin{pmatrix} I_n &&Z \\ &I_n& \\ & &1 \end{pmatrix}  \begin{pmatrix} {^tg^{-1}} && \\ &{^tg^{-1}}& \\ & &1 \end{pmatrix}  \tau_{2n+1} \right) & \psi(\mathrm{Tr}(X))\\
 & \times dX dZ|\mathrm{det}(g)|.
\end{split}
\]
We insert the definition of $\varphi_1$ and $\varphi_2$ and then utilize the assignment $g \mapsto \tau_{2n+1}{^{\iota}g}\tau^{-1}_{2n+1}$ on the last matrix. After the change of the variables $y \mapsto -y$, the identity becomes
\[
\begin{split}
 &\int_{\mathbb{A}^n}\int_{k^n \backslash \mathbb{A}^n} \int_{\mathcal{M}_n(k) \backslash \mathcal{M}_n(\mathbb{A}) }  \varphi \left(  \begin{pmatrix} I_n &X&Z \\ &I_n& \\ & &1 \end{pmatrix}
   \begin{pmatrix} g && \\ &g& \\ & &1 \end{pmatrix} \begin{pmatrix} I_n && \\ &I_n& \\ &y &1 \end{pmatrix}
\right) \psi^{-1}(\mathrm{Tr}(X))\Phi(y) dX dZ dy\\
&=\int_{\mathbb{A}^n}\int_{k^n \backslash \mathbb{A}^n} \int_{\mathcal{M}_n(k) \backslash \mathcal{M}_n(\mathbb{A}) }  \widetilde{\varphi} \left( 
\begin{pmatrix} I_n &X&Z \\ &I_n& \\ & &1 \end{pmatrix}
 \begin{pmatrix} {^tg^{-1}} && \\ &{^tg^{-1}}& \\ & &1 \end{pmatrix}
  \begin{pmatrix} I_n && \\ &I_n& \\ &y &1 \end{pmatrix}
 \tau_{2n+1} \right) \\
& \quad \times \psi(\mathrm{Tr}(X))\hat{\Phi}(y) dX dZ dy |\mathrm{det}(g)|
\end{split}
\]
from which the desired global functional equation for integrals follows.
\end{proof}

Let
\[
  W_{\varphi}(g)=\int_{N_m(k) \backslash N_m(\mathbb{A}) } \varphi(ng) \psi^{-1}(n) dn \quad \text{and} \quad
  \widetilde{W}_{\varphi}(g)=\int_{N_m(k) \backslash N_m(\mathbb{A}) } \widetilde{\varphi}(w_mng)  \psi(n) dn
\]
be the associated Whittaker function of $\varphi$ and $\widetilde{\varphi}$, respectively. We have yet to check that our integrals are Eulerian.

\begin{proposition}[Jacquet-Shalika]
\label{J-Sunfolding}
For $\varphi \in V_{\Pi}$ and $\Phi \in \mathcal{S}(F^n)$, global Jacquet-Shalika integrals 
\[
\begin{split}
  &J_{\psi}(s,W_{\varphi},\Phi)\\
  &=\int_{N_n(\mathbb{A}) \backslash GL_n(\mathbb{A})} \int_{\mathcal{N}_n(\mathbb{A}) \backslash \mathcal{M}_n(\mathbb{A})}
  W_{\varphi} \left( \begin{pmatrix} I_n & X \\ & I_n \end{pmatrix} \begin{pmatrix} g &  \\ & g \end{pmatrix} \right)  \psi^{-1}({\rm Tr}(X)) \Phi(e_ng) |\mathrm{det}(g)|^s dX dg
\end{split}
\]
in the even case $m=2n$ and
\[
\begin{split}
  J_{\psi}(s,W_{\varphi},\Phi)
  &=\int_{N_n(\mathbb{A}) \backslash GL_n(\mathbb{A})} \int_{\mathcal{N}_n(\mathbb{A}) \backslash \mathcal{M}_n(\mathbb{A})} \int_{\mathbb{A}^n}
  W_{\varphi} \left( \begin{pmatrix} I_n & X &  \\ & I_n &\\ &&1 \end{pmatrix} \begin{pmatrix} g &&  \\ & g& \\ &&1 \end{pmatrix} \begin{pmatrix} I_n &  &  \\ & I_n &\\ &y&1 \end{pmatrix}  \right) \\
  &\quad \times \psi^{-1}({\rm Tr}(X)) \Phi(y) |\mathrm{det}(g)|^{s-1} dy dX dg
\end{split}
\]
in the odd case $m=2n+1$ converges when $\mathrm{Re}(s)$ is sufficiently large and, when this is the case, we have
\[
  I_{\psi}(s,\varphi,\Phi)= J_{\psi}(s,W_{\varphi},\Phi).
\]
We suppose, in addition, that $W_{\varphi}(g)=\prod_{v \in |\mathbb{P}^1|}W_{\varphi_v}(g_v)$, $\psi(n)=\prod_{v \in |\mathbb{P}^1|} \psi(n_v)$, and $\Phi(g)=\prod_{v \in |\mathbb{P}^1|} \Phi_v(g_v)$. Then, when $\mathrm{Re}(s)$ is sufficiently large,
\[
 J_{\psi}(s,W_{\varphi},\Phi)=\prod_{v \in |\mathbb{P}^1|}  J_{\psi_v}(s,W_{\varphi_v},\Phi_v).
\]
Likewise, the right hand side of the functional equation is also unfold and can be factored as
\[
 I_{\psi^{-1}}(1-s,\varrho(\tau_m)\widetilde{\varphi},\hat{\Phi})=J_{\psi^{-1}}(1-s,\varrho(\tau_m)\widetilde{W}_{\varphi},\hat{\Phi})
 =\prod_{v \in |\mathbb{P}^1|} J_{\psi^{-1}}(1-s,\varrho(\tau_m)\widetilde{W}_{\varphi_v},\hat{\Phi}_v)
\]
with the convergence for $\mathrm{Re}(s) \ll 0$.
\end{proposition}

\begin{proof}
All these statements are drawn, with some minor changes of notation, from Proposition 5 in \S 6 of \cite{JaSh88} for the even case $m=2n$, and Section 9.2 of \cite{JaSh88} 
for the odd case $m=2n+1$.
\end{proof}

Throughout, we will take $S \subset |\mathbb{P}^1|$ to be a finite set of places such that for all $v \notin S$,
$\Pi_v$ and $\psi_v$ are all unramified and $\psi_v$ normalized. The partial $L$-function is a product of local factors
\[
  L^S(s,\Pi,\wedge^2)=\prod_{v \notin S} L(s,\Pi_v,\wedge^2).
\]
More precisely, this product converges for $\mathrm{Re}(s)$ large enough (cf.\cite[\S 8, \S 9]{JaSh88}). The
global $L$-function and $\varepsilon$-factors for $\Pi$ are
\[
 L(s,\Pi,\wedge^2,S)=\prod_{v \in |\mathbb{P}^1|} L(s,\Pi_v,\wedge^2)=L^S(s,\Pi,\wedge^2)\prod_{v \in S} L(s,\Pi_v,\wedge^2)
\]
and
\[
\varepsilon(s,\Pi,\wedge^2,S)=\prod_{v \in |\mathbb{P}^1|} \varepsilon(s,\Pi_v,\wedge^2,\psi_v)=\prod_{v \in S} \varepsilon(s,\Pi_v,\wedge^2,\psi_v).
\]
As for the $\varepsilon$-factor, we know that $\varepsilon(s,\Pi_v,\wedge^2,\psi_v) \equiv 1$ for $v \notin S$. The independence of $\varepsilon(s,\Pi,\wedge^2,S)$ from the choice of $\psi$ can be seen as a consequence of the global functional equation below.

\begin{theorem}
\label{global-func-RS}
The global $L$-function $L(s,\Pi,\wedge^2,S)$ has a meromorphic continuation to the entire plane and it satisfies the global functional equation 
\[
  L(s,\Pi,\wedge^2,S)=\varepsilon(s,\Pi,\wedge^2,S)L(1-s,\widetilde{\Pi},\wedge^2,S),
\]
where $\varepsilon(s,\Pi,\wedge^2,S)$ is entire and non-vanishing. This identity implies that $\varepsilon(s,\Pi,\wedge^2,S)$ is independent of $\psi$ as well. 
\end{theorem}

\begin{proof}
From the unfolding in Proposition \ref{J-Sunfolding}, and the local calculation of \cite[\S 7.2, \S 9.4]{JaSh88} together with Corollary \ref{unramified}, we know that for  
$\mathrm{Re}(s)$ large and for appropriate choice of $\varphi$, we have
\[
\begin{split}
  I_{\psi}(s,\varphi,\Phi)&= J_{\psi}(s,W_{\varphi},\Phi)=\prod_{v \in |\mathbb{P}^1|}  J_{\psi_v}(s,W_{\varphi_v},\Phi_v)
  =\left( \prod_{v \in S} J_{\psi_v}(s,W_{\varphi_v},\Phi_v) \right) L^S(s,\Pi,\wedge^2)\\
  &=\left( \prod_{v \in S} \frac{J_{\psi_v}(s,W_{\varphi_v},\Phi_v)}{L(s,\Pi_v,\wedge^2)} \right) L(s,\Pi,\wedge^2,S)
  =\left(  \prod_{v \in S}  e_v(s,W_{\varphi_v},\Phi_v) \right)L(s,\Pi,\wedge^2,S),
\end{split}
\]
where
$
 e_v(s,W_{\varphi_v},\Phi_v)=\dfrac{J_{\psi_v}(s,W_{\varphi_v},\Phi_v)}{L(s,\Pi_v,\wedge^2)}.
$
It can be seen from Theorem \ref{exterior-structure} that $e_v(s,W_{\varphi_v},\Phi_v)$ is entire. Therefore $L(s,\Pi,\wedge^2,S)$ has a meromorphic continuation as the integral 
$I_{\psi}(s,\varphi,\Phi)$ is a meromorphic function on the entire plane from Theorem \ref{global-func}.
While on the other side, we obtain
\[
\begin{split}
   I_{\psi^{-1}}(1-s,\varrho(\tau_m)\widetilde{\varphi},\hat{\Phi})&=J_{\psi^{-1}}(1-s,\varrho(\tau_m)\widetilde{W}_{\varphi},\hat{\Phi})\\
   &=\left( \prod_{v \in S} \widetilde{e}_v(1-s,\varrho(\tau_m)\widetilde{W}_{\varphi_v},\hat{\Phi}_v)  \right)L(1-s,\widetilde{\Pi},\wedge^2,S)
\end{split}
\]
with $\widetilde{e}_v(1-s,\varrho(\tau_m)\widetilde{W}_{\varphi_v},\hat{\Phi}_v)=\dfrac{J_{\psi^{-1}}(1-s,\varrho(\tau_m)\widetilde{W}_{\varphi_v},\hat{\Phi}_v)}{L(1-s,\widetilde{\Pi}_v,\wedge^2)}$. However we derive from the local functional equation, Theorem \ref{exterior-func}, that
\[
\begin{split}
 \widetilde{e}_v(1-s,\varrho(\tau_m)\widetilde{W}_{\varphi_v},\hat{\Phi}_v)
 &=\frac{J_{\psi^{-1}}(1-s,\varrho(\tau_m)\widetilde{W}_{\varphi_v},\hat{\Phi}_v)}{L(1-s,\widetilde{\Pi}_v,\wedge^2)}
 =\varepsilon(s,\Pi_v,\wedge^2,\psi_v)\frac{J_{\psi_v}(s,W_{\varphi_v},\Phi_v)}{L(s,\Pi_v,\wedge^2)}\\
 &=\varepsilon(s,\Pi_v,\wedge^2,\psi_v)e_v(s,W_{\varphi_v},\Phi_v).
\end{split}
\]
Combining these all together, we get
\[
 L(s,\Pi,\wedge^2,S)=\left( \prod_{v \in S} \varepsilon(s,\Pi_v,\wedge^2,\psi_v) \right) L(1-s,\widetilde{\Pi},\wedge^2,S)=\varepsilon(s,\Pi,\wedge^2,S)L(1-s,\widetilde{\Pi},\wedge^2,S)
\]
since for $v \notin S$ we know that $\Pi_v$ and $\psi_v$ are unramified so that $\varepsilon(s,\Pi_v,\wedge^2,\psi_v) \equiv 1$.
\end{proof}

\subsection{The equality for discrete series representations}
\label{sec3.2}

Let $k_0=\mathbb{F}_q((t))$ be the completion of $k$ at the point $0 \in |\mathbb{P}^1|$. We start with a local irreducible unitary supercuspidal 
representation $\rho_{\circ}$ and globalize it according to the result of Henniart-Lomel\'{\i}  \citelist{\cite{HL-Exterior}\cite{HL-RS}*{Theorem 3.1}}.

\begin{theorem}[Henniart-Lomel\'{\i}]
\label{HL-global}
 Let $\rho_{\circ}$ be an irreducible unitary supercuspidal representation of $GL_m(F)$. We choose an isomorphism $\xi : F \overset{\sim}{\rightarrow} k_0$. Then there exists a cuspidal unitary automorphic representation
$\Pi=\bigotimes'_v \Pi_v$ whose local components $\Pi_v$ satisfy:
\begin{itemize}
\item $\rho_{\circ}$ corresponds to $\Pi_0$ via $\xi$;
\item at the places $v \in |\mathbb{P}^1|$ away from $0$, $1$, and $\infty$, $\Pi_v$ is irreducible and unramified;
\item $\Pi_1$ is an irreducible generic subquotient of an unramified principal series representation;
\item $\Pi_{\infty}$ is an irreducible generic subquotient of a tamely ramified principal series representation.
\end{itemize}
\end{theorem}

We have control at all places outside $0$, which makes it possible to deduce the identity for irreducible supercuspidal representations.

\begin{theorem}[Supercuspidal cases]
\label{self-dual-supercusp} 
Let $\rho$ be an irreuducible supercuspidal representation of $GL_r$. Then we have
\[
  L(s,\rho,\wedge^2)=\mathcal{L}(s,\rho,\wedge^2).
\]
As a consequence, if $\rho$ is $(S_{2n},\Theta)$-distinguished, then $\rho$ is self-dual.
\end{theorem}

\begin{proof} 
Twisting by an unramified character does not affect the conclusion, so we can assume that $\rho=\rho_{\circ}$ is unitary. (See the proof of Theorem \ref{BF-disc} for details, cf.\cite[\S 6.6]{Lomeli}).
We define the Langlands-Shahidi global $L$-function and $\varepsilon$-factors for $\Pi$ by
\[
  \mathcal{L}(s,\Pi,\wedge^2,S)=\prod_{v \in |\mathbb{P}^1|} \mathcal{L}(s,\Pi_v,\wedge^2) \;\; \text{and} \;\;
  \varepsilon_{LS}(s,\Pi,\wedge^2,\psi,S)=\prod_{v \in |\mathbb{P}^1|} \varepsilon_{LS}(s,\Pi_v,\wedge^2,\psi_v).
\]
We choose a finite set $S=\{0,1,\infty \}$ of places. Applying Theorem \ref{HL-global} to the irreducible unitary supercuspidal representation $\rho_{\circ}$, we obtain a cuspidal unitary automorphic representation $\Pi$. For our convenience, we rewrite the global functional equation in \cite[4.1.(vi)]{HL-Exterior} as
\begin{equation}
\label{LS-func}
  \mathcal{L}(s,\Pi,\wedge^2,S)=\varepsilon_{LS}(s,\Pi,\wedge^2,S)\mathcal{L}(1-s,\widetilde{\Pi},\wedge^2,S).
\end{equation}
The function $\varepsilon_{LS}(s,\Pi,\wedge^2,S)$ is entire and non-vanishing. From the global functional equation
given by Theorem \ref{global-func-RS} and \eqref{LS-func}, this means that the ratio of $L$-function satisfies
\[
 \frac{L(s,\Pi,\wedge^2,S)}{\mathcal{L}(s,\Pi,\wedge^2,S)}=\eta(s,\Pi,S)\frac{L(1-s,\widetilde{\Pi},\wedge^2,S)}{\mathcal{L}(1-s,\widetilde{\Pi},\wedge^2,S)},
\]
where $\eta(s,\Pi,S)=\varepsilon(s,\Pi,\wedge^2,S)\varepsilon_{LS}(s,\Pi,\wedge^2,S)^{-1}$ is entire and non-vanishing. Applying the already established 
principal series representations in Corollary \ref{unramified} alongside with Proposition \ref{principal} at the places $|\mathbb{P}^1|-\{0\}$ yields agreements:
\[
 \prod_{v \notin \{0\}} L(s,\Pi_v,\wedge^2)=
  \prod_{v \notin \{0\}} \mathcal{L}(s,\Pi_v,\wedge^2)
\quad
\text{and}
\quad
  \prod_{v \notin \{0\}} L(1-s,\widetilde{\Pi}_v,\wedge^2)=
  \prod_{v \notin \{0\}} \mathcal{L}(1-s,\widetilde{\Pi}_v,\wedge^2).
\]
Therefore, at the remaining one place, we have
\[
  \frac{ L(s,\rho_{\circ},\wedge^2)}{\mathcal{L}(s,\rho_{\circ},\wedge^2)}=\eta(s,\Pi,S)\frac{L(1-s,\widetilde{\rho}_{\circ},\wedge^2)}{\mathcal{L}(1-s,\widetilde{\rho}_{\circ},\wedge^2)}.
\]
In view of \cite{GL} and \cite[Theorem 3.7]{KeRa}, $L(s,\rho_{\circ},\wedge^2)$ and $\mathcal{L}(s,\rho_{\circ},\wedge^2)$ are regular and non-vanishing in the region $\mathrm{Re}(s) > 0$,
whereas similar analytic properties for $L(1-s,\widetilde{\rho}_{\circ},\wedge^2)$ and $\mathcal{L}(1-s,\widetilde{\rho}_{\circ},\wedge^2)$ are valid in the half plane $\mathrm{Re}(s) < 1$. 
This forces that the ratio $ \dfrac{ L(s,\rho_{\circ},\wedge^2)}{\mathcal{L}(s,\rho_{\circ},\wedge^2)}$ is an entire and non-vanishing function and hence it is a constant. Since these $L$-factors
are normalized, these must be equal.

\par
 We now gain the full strength of flexibility to transport $L$-factors in the Langlands-Shahidi side to the Rankin-Selberg side. The $L$-factor $\mathcal{L}(s,\rho \times \rho)$ is decomposed as the product of $\mathcal{L}(s,\rho,\wedge^2)$ and $\mathcal{L}(s,\rho,\mathrm{Sym}^2)$ (cf. \citelist{\cite{GL}\cite{HL-Exterior}\cite{Shahidi}*{Corollary 8.2}}). Then the pole of $\mathcal{L}(s,\rho,\wedge^2)$ at $s=0$ detected by the existence of Shalika functional \cite[Theorem 3.6.(ii)]{Jo20-2} contributes the pole of $\mathcal{L}(s,\rho \times \rho)$. This is amount to saying that $\rho$ is self-dual \cite[Proposition 4.6]{Matringe15}.
\end{proof}

Once we know the inductivity of $\varepsilon$-factors, we expect that $\eta(s,\Pi,S) \equiv 1$ independent the choice of $S$. 
We now come to the case of discrete series representations.

\begin{theorem}[Discrete series cases] 
\label{disc-ext-selfdual}
Let $\Delta$ be a discrete series representation of $GL_m$. Then we have
\[
  L(s,\Delta,\wedge^2)=\mathcal{L}(s,\Delta,\wedge^2).
\]
As a consequence, if $\Delta$ is $(S_{2n},\Theta)$-distinguished, then $\Delta$ is self-dual.
\end{theorem}

\begin{proof}
As indicated in the proof of Theorem \ref{BF-disc}, after proper unramified twisting of $\Delta$, we can easily reduce the equality to the case when $\Delta$ is a unitary representation of the form $[\rho_{\circ}\nu^{-\frac{\ell-1}{2}}, \dotsm,  \rho_{\circ}\nu^{\frac{\ell-1}{2}}]$ with $\rho_{\circ}$  a unitary irreducible supercuspidal representation of $GL_r$ (cf.\cite[\S 6.6]{Lomeli}). 
Taking the advantage of Proposition \ref{self-dual-supercusp}, Theorem \ref{Discrete} finally matches with the expression in \cite[Proposition 8.1]{Shahidi}.
Concerning the second assertion, we literarily reiterate the second part of the proof of Theorem \ref{self-dual-supercusp} line-by-line, and therefore we omit thorough arguments entirely. 
\end{proof}



The identity can be extended to the class of all irreducible admissible representations of  $GL_m$. 

\begin{theorem}
Let $\pi$ be an irreducible admissible representation of $GL_m$. Then we have
\[
  L(s,\pi,\wedge^2)=\mathcal{L}(s,\pi,\wedge^2).
\]
\end{theorem}

\begin{proof}
We realize $\pi$ as the unique Langlands quotient of Langlands type $\Xi=\mathrm{Ind}^{GL_m}_{\rm Q}({\Delta_{\circ}}_1\nu^{u_1} \otimes {\Delta_{\circ}}_2\nu^{u_2} \otimes \dotsm \otimes {\Delta_{\circ}}_t\nu^{u_t})$ which is again of Whittaker type. 
Thanks to Theorem \ref{disc-ext-selfdual}, the Working Hypothesis is not required to be checked for discrete series representations. Then Theorem \ref{ext-langlands-prod} gives us that 
 \[
  L(s,\Xi,\wedge^2)=\prod_{1 \leq k \leq t} L(s+2u_k,{\Delta_{\circ}}_k,\wedge^2) \prod_{1 \leq i < j  \leq t} L(s+u_i+u_j, {\Delta_{\circ}}_i \times {\Delta_{\circ}}_j),
\]
 which coincides with corresponding decompositions in the Langlands-Shahidi theory \cite[\S 3.1.(xi)]{GL}.
\end{proof}

By exploiting the main result of Henniart-Lomel\'{\i} \cite{HL-Exterior}, it can be summarized that the definition of local exterior square $L$-function via the theory of integral representations is compatible with the local Langlands correspondence. In what follows, we let $W'_F$ denote the Weil–Deligne group of $F$ and let $\phi$ a $m$-dimensional (complex-valued) Frobenius semi-simple representation of $W'_F$. We call this 
the {\it Weil-Deligne representation} of $W'_F$. Let $\wedge^2$ denote the exterior representation of $GL_m(\mathbb{C})$. We then denote by $L(s,\wedge^2 ( \phi))$ the {\it Artin exterior square} $L$-{\it factor} attached to $\phi$. 

\begin{theorem}
\label{exter-equal}
Let $\pi$ be an irreducible admissible representation of $GL_m(F)$ and $\phi(\pi)$ the Weil-Deligne representation of $W'_F$ corresponding to $\pi$ under the local Langlands correspondence. Then we have
\[
  L(s,\pi,\wedge^2)=\mathcal{L}(s,\pi,\wedge^2)=L(s,\wedge^2( \phi(\pi))).
\]
\end{theorem}

\section{Bump-Friedberg and Flicker Zeta Integrals}

\subsection{Bump-Friedberg $L$-factors} 
\label{sec:BF}
We define the embedding $J: GL_n \times GL_{n} \rightarrow GL_m$ by
\[
 J(g,g^{\prime})_{k,l}=
  \begin{cases}
  g_{i,j} &  \text{if $k=2i-1$, $l=2j-1$, } \\
  g^{\prime}_{i,j} &  \text{if $k=2i$, $l=2j$,} \\
 0 &  \text{otherwise,} \\
  \end{cases}
\]
for $m=2n$ even
and  $J: GL_{n+1} \times GL_n \rightarrow GL_m$ by
\[
 J(g,g^{\prime})_{k,l}=
  \begin{cases}
  g_{i,j} &  \text{if $k=2i-1$, $l=2j-1$, } \\
  g^{\prime}_{i,j} &  \text{if $k=2i$, $l=2j$,} \\
 0 &  \text{otherwise,} \\
  \end{cases}
\]
for $m=2n+1$ odd. As for the intention of holding onto coherent terminology with \cite{Matringe15}, 
interested readers may perceive that we interchange the role of $g$ and $g'$ in \cite{BF}. Let $\pi=\mathrm{Ind}^{GL_m}_{\rm Q}(\Delta_1 \otimes \dotsm \otimes \Delta_t)$ be a parabolically induced representation. For each Whittaker function $W \in \mathcal{W}(\pi,\psi)$ and Schwartz-Bruhat function $\Phi \in \mathcal{S}(F^n)$, we define Bump-Friedberg integrals:
\[
 Z(s_1,s_2,W,\Phi)=\int_{N_{n} \backslash GL_{n}} \int_{N_n \backslash GL_n} W(J(g,g^{\prime})) \Phi(e_mJ(g,g^{\prime})) |\mathrm{det}(g)|^{s_1-1/2} |\mathrm{det}(g^{\prime})|^{1/2+s_2-s_1} dg dg^{\prime}
\]
in the even case $m=2n$ and
\[
 Z(s_1,s_2,W,\Phi)=\int_{N_{n} \backslash GL_{n}} \int_{N_{n+1} \backslash GL_{n+1}} W(J(g,g^{\prime})) \Phi(e_mJ(g,g^{\prime})) |\mathrm{det}(g)|^{s_1} |\mathrm{det}(g^{\prime})|^{s_2-s_1} dg dg^{\prime}
 \]
 in the odd case $m=2n+1$. If $r$ a real number, we denote by $\delta_r$ the character,
 \[
  \delta_r : J(g,g') \mapsto   \left| \frac{\mathrm{det}(g)}{\mathrm{det}(g')} \right|^r.
 \]
  We denote by $\chi_m$ and $\mu_m$ characters of $H_m$:
\[
  \chi_m\left( w_m \begin{pmatrix} g & \\ & g' \end{pmatrix} w^{-1}_m\right)=
  \begin{cases}
{\bf 1}_{H_m} & \text{for $m=2n$;} \\
\displaystyle  \left| \frac{\mathrm{det}(g)}{\mathrm{det}(g^{\prime})} \right| & \text{for $m=2n+1$;} \\
\end{cases}
\]
\[
 \mu_m\left( w_m \begin{pmatrix} g & \\ & g' \end{pmatrix} w^{-1}_m\right)=
  \begin{cases}
\displaystyle  \left| \frac{\mathrm{det}(g)}{\mathrm{det}(g^{\prime})} \right| & \text{for $m=2n$;} \\
{\bf 1}_{H_m} & \text{for $m=2n+1$.} \\
\end{cases}
\]
We turn toward the case for $s_1=s+1/2$ and $s_2=2s$. We unify Bump-Friedberg Zeta Integrals as one single integral of the form:
\[
Z(s,W,\Phi)=\int_{(N_m \cap H_m) \backslash H_m} W(h)  \chi^{1/2}_m(h) \Phi(e_mh) |\mathrm{det}(h)|^s \, dh
\]
The twisted analogue of Bump-Friedberg Zeta Integrals attached to $\chi_{\alpha}$ is defined by
\[
Z(s,W,\Phi,\chi_{\alpha})=\int_{(N_m \cap H_m) \backslash H_m} W(h) \chi_{\alpha}(h) \chi^{1/2}_m(h) \Phi(e_mh) |\mathrm{det}(h)|^s \, dh
\]
The integral $Z(s,W,\Phi,\chi_{\alpha})$ converges absolutely for $s$ of real part large enough. The $\mathbb{C}$-vector space generated by
Bump-Friedberg Zeta integrals
\[ 
  \langle Z(s,W,\Phi,\chi_{\alpha}) \; | \;  W \in \mathcal{W}(\pi,\psi), \Phi \in \mathcal{S}(F^m) \rangle
\]
is a $\mathbb{C}[q^{\pm s}]$-fractional ideal $\mathcal{I}(\pi,\chi_{\alpha},BF)$ of $\mathbb{C}(q^{-s})$. The ideal $\mathcal{I}(\pi,\chi_{\alpha},BF)$ is principal 
and has a unique generator of the form $P(q^{-s})^{-1}$, where $P(X)$ is a polynomial in $\mathbb{C}[X]$ with $P(0)=1$.
The {\it Bump-Friedberg} $L$-{\it factor} associated to $\pi$ is defined by the unique normalized generator \cite[Proposition 4.8]{Matringe15}:
\[
  L(s,\pi,\chi_{\alpha},BF)=\frac{1}{P(q^{-s})}.
\]
If $\alpha=\textbf{1}_{F^{\times}}$ is a trivial character, it is convenient to write $L(s,\pi,BF)$ for $L(s,\pi,\chi_{\textbf{1}_{F^{\times}}},BF)$. The {\it Bump-Friedberg} $\gamma$-{\it factor}
\[
  \gamma(s,\pi,BF,\psi) =\varepsilon(s,\pi,BF,\psi) \frac{L(1/2-s,\pi^{\iota},\delta_{-1/2},BF)}{L(s,\pi,BF)}
\]
is a rational function in $\mathbb{C}(q^{-s})$ that depends on a choice of a non-trivial character $\psi$ (Cf. \cite[Proposition 4.11]{Matringe15}). 
While the proof of \cite[Proposition 6.2]{Matringe} reflects the structure of Weil-Deligne representations, our aim is to show the factorization of 
$L(s,\pi,\chi_{\alpha},BF)$ as a product of the standard $L$-factor $L(s+1/2,\pi)$ and the exterior square $L$-factor $L(2s,\pi,\wedge^2)$
within the framework of Rankin-Selberg method. Our approach here is more direct and concise.

\begin{theorem}[Supercuspidal cases] 
\label{BF-supercusp-factor}
Let $\rho$ be an irreducible supercuspidal representation of $GL_r$. Then
\[
 L(s,\rho,BF)=L(s+1/2,\rho)L(2s,\rho,\wedge^2).
\]
\end{theorem}

\begin{proof}
If $r=1$, then $\rho$ is a character of $F^{\times}$. The integral is just Tate integrals of the form
$\displaystyle
 \int_{F^{\times}} \rho(z) \Phi(z) |z|^{s+\frac{1}{2}} d^{\times} z,
$
hence 
\[
L(s,\rho,BF)=L(s+1/2,\rho)=L(s+1/2,\rho)L(2s,\rho,\wedge^2),
\]
 where the last equality comes from $L(2s,\rho,\wedge^2)=1$ (Cf. \cite[Theorem 2.13]{Jo20-2}).

\par
We deduce from Theorem \ref{factor} aligned with \cite[Proposition 4.14]{Matringe15} that all the poles of $L(s,\rho,BF)$
and $L(s,\rho,\wedge^2)$ are necessarily simple. Given $r=2n+1$ with $n \geq 1$, the result of Matringe \cite[Theorem 3.1]{Matringe} (Cf. Theorem \ref{M-Shalika-odd})
tells us that $\rho$ cannot be $H_{2n+1}$-distinguished. According to \cite[\S 3.1]{Jacquet} coupled with \cite[Theorem 3.6.(ii)]{Jo20-2} and \cite[Corollary 4.3]{Matringe15}, 
we have
\[
  L(s,\rho,BF)=L(s,\rho)=L(s,\rho,\wedge^2)=1.
\]

\par
Now we turn to the case when $r=2n$ is even. Analyzing poles of local $L$-functions is just a question of certain distinctions of representations. To be precise, \cite[Theorem 3.6.(i)]{Jo20-2} together with \cite[Corollary 4.3]{Matringe15} and \S \ref{sec:2.2} lead us to the following equivalent statements:
\begin{enumerate}[label=$(\mathrm{\arabic*})$]
\item\label{Char-pole-1} $L(2s,\rho,\wedge^2)$ has a pole at $s=s_0$;
\item\label{Char-pole-2} $L(s,\rho,BF)$ has a pole at $s=s_0$;
\item\label{Char-pole-3} $\rho\nu^{s_0}$ is $(S_{2n},\Theta)$-distinguished;
\item\label{Char-pole-3} $\rho\nu^{s_0}$ is $H_{2n}$-distinguished.
\end{enumerate} 
The above characterization of poles of $L$-factors can be reinterpreted in the following manner:
\[
  L(s,\rho,BF)=L(2s,\rho,\wedge^2)=L(s+1/2,\rho)L(2s,\rho,\wedge^2),
\]
where the last identity follows from $L(s+1/2,\rho)=1$ (Cf. \cite[\S 3.1]{Jacquet}).
\end{proof}

Unlike the case of Jacquet and Shalika's Zeta integrals \S \ref{sec2}-\S \ref{sec3}, it is necessary to use the hereditary property of $H_{2m+1}$-distinguished representations due to N. Matringe \cite[Theorem 3.1]{Matringe15} additionally.

\begin{theorem}[N. Matringe, Odd cases $m=2n+1$]
\label{M-Shalika-odd}
Let $\pi=\mathrm{Ind}^{GL_{2n+1}}_{\rm Q}(\Delta_1 \otimes \Delta_2 \otimes \dotsm \otimes \Delta_t)$ be an irreducible generic representation of $GL_{2n+1}$.
Let $\alpha$ be a character of $F^{\times}$ with $0 \leq \mathrm{Re}(\alpha) \leq 1/2$. Then $\pi$ is $(H_{2n+1},\chi_{\alpha}\delta_{-1/2})$-distinguished if and \
only if $\pi$ is a parabolically induced representation of the form $\mathrm{Ind}^{GL_{2n+1}}_{P_{2n,1}}(\pi' \otimes \alpha\nu^{-1/2})$, for $\pi'$ an irreducible generic
$(H_{2n},\chi_{\alpha})$-distinguished representation of $GL_{2n}$ such that $\mathrm{Ind}^{GL_{2n+1}}_{P_{2n,1}}(\pi' \otimes \alpha\nu^{-1/2})$ is still irreducible and generic.
\end{theorem}

 Throughout the rest of \S \ref{sec:BF}, a variant of the systematic machinery developed in \S \ref{deform-special} (in particular, Proposition \ref{three} to Theorem \ref{Deform-gamma}) should continue to work out 
in the context of Bump-Friedberg Zeta integrals and it is dealt with in \cite[\S 4]{Matringe15} in great detail and clarity. By doing so, Bump-Friedberg local $L$-functions are compatible with
the classification of discrete series representation in terms of supercuspidal ones owing to Bernstein and Zelevinsky \cite{BeZe,Zelevinsky}.

\begin{proposition}
\label{Discrete-BF}
Let $\Delta_{\circ}=[\rho_{\circ}\nu^{-\frac{\ell-1}{2}},\dotsm,\rho_{\circ}\nu^{\frac{\ell-1}{2}}]$ be an irreducible square integrable representation of $GL_{\ell r}$
with $\rho_{\circ}$ an irreducible unitary supercuspidal representation of $GL_r$.
 \begin{enumerate}[label=$(\mathrm{\roman*})$]
\item\label{Discrete-BF-1} Suppose that $\ell$ is even. Then we have
\[
 L(s,\Delta_{\circ},BF)=\displaystyle
L(s+\ell/2,\rho_{\circ} )\prod_{i=1}^{\ell/2}L(2s,\rho_{\circ}\nu^{(\ell+1)/2-i},\wedge^2)L(2s,\rho_{\circ}\nu^{\ell/2-i},\mathrm{Sym}^2).
\]
\item\label{Discrete-BF-2} Suppose that $\ell$ is odd. Then we have
\[
L(s,\Delta_{\circ},BF)=\displaystyle
L(s+\ell/2,\rho_{\circ} )\prod_{i=1}^{(\ell+1)/2}L(2s,\rho_{\circ}\nu^{(\ell+1)/2-i},\wedge^2)\prod_{i=1}^{(\ell-1)/2}L(2s,\rho_{\circ}\nu^{\ell/2-i},\mathrm{Sym}^2). 
\]
 \end{enumerate}
\end{proposition}

\begin{proof}
By the uniqueness of the Whittaker functional, the Whittaker model for $\Delta_{\circ}$
coincides with that for $\mathrm{Ind}(\rho_{\circ}\nu^{-\frac{\ell-1}{2}}\otimes \dotsm \otimes \rho_{\circ}\nu^{\frac{\ell-1}{2}})$.
Likewise the same trait holds for dual objects provided by $\mathrm{Ind}(\widetilde{\rho}_{\circ}\nu^{-\frac{\ell-1}{2}}\otimes \dotsm \otimes \widetilde{\rho}_{\circ}\nu^{\frac{\ell-1}{2}})$ and $\widetilde{\Delta}_{\circ}$. According to \cite[Proposition 5.5]{Matringe15},
\[
\gamma(s,\Delta_{\circ},BF,\psi)
 \sim \prod_{i=0}^{\ell-1} \gamma \left( s+\frac{1-\ell}{2}+i,\rho_{\circ},BF,\psi \right) \prod_{0 \leq i < j \leq \ell-1} \gamma(2s+1-\ell+i+j,\rho_{\circ} \times \rho_{\circ}, \psi).
\]
With the help of Theorem \ref{BF-supercusp-factor}, the expression can be reformulated in terms of $L$-factors as
\[
\begin{split}
 &\gamma(s,\Delta_{\circ},BF,\psi) \\
 & \sim  \prod_{i=0}^{\ell-1} \frac{L(-s-i+\ell/2,\widetilde{\rho}_{\circ})}{L(s+i+1-\ell/2,\rho_{\circ})}  
  \prod_{i=0}^{\ell-1} \frac{L(-2s+\ell-2i,\widetilde{\rho}_{\circ},\wedge^2)}{L(2s+1-\ell+2i,\rho_{\circ},\wedge^2)}
\prod_{0 \leq i < j \leq \ell-1} \frac{L(-2s+\ell-i-j,\widetilde{\rho}_{\circ} \times \widetilde{\rho}_{\circ})}{L(2s+1-\ell+i+j,\rho_{\circ} \times \rho_{\circ})}.
\end{split}
\]
By virtue of \cite[Lemma 5.11]{Jo20-2} combined with $L(-s,\rho_{\circ}) \sim L(s, \widetilde{\rho}_{\circ})$ (Cf. \cite[\S 8.2 (15), (16)]{JPSS}), it may be written as
\[
\begin{split}
 & \gamma(s,\Delta_{\circ},BF,\psi) \\
 & \sim  \prod_{i=0}^{\ell-1} \frac{L(s+i-\ell/2,\rho_{\circ})}{L(s+i+1-\ell/2,\rho_{\circ})}  
 \prod_{i=0}^{\ell-1} \frac{L(2s-\ell+2i,\rho_{\circ},\wedge^2)}{L(2s+1-\ell+2i,\rho_{\circ},\wedge^2)}
\prod_{0 \leq i < j \leq \ell-1} \frac{L(2s-\ell+i+j,\rho_{\circ} \times \rho_{\circ})}{L(2s+1-\ell+i+j,\rho_{\circ} \times \rho_{\circ})}.
\end{split}
\]
We do the case $\ell$ even, the case $\ell$ odd being treated similarly. At this moment, we repeat the proof given in \cite[Theorem 5.12]{Jo20-2} with adjusting $s$ to $2s$. After cancelling common factors, our quotient is simplified to
\begin{equation}
\label{reduce}
\gamma(s,\Delta_{\circ},BF,\psi) \sim \frac{L(s-\ell/2,\rho_{\circ})}{L(s+\ell/2,\rho_{\circ})} \prod_{i=0}^{(\ell/2)-1} \frac{L(2s-\ell+2i,\rho_{\circ},\wedge^2)L(2s-\ell+2i+1,\rho_{\circ},\mathrm{Sym}^2)}{L(2s+2i+1,\rho_{\circ},\wedge^2)L(2s+2i,\rho_{\circ},\mathrm{Sym}^2)}. 
\end{equation}
Using \cite[Corollary 4.1]{Matringe15}, $L(1/2-s,\widetilde{\Delta}_{\circ},\delta_{-1/2},BF)^{-1}$ has zeros in the half plane $\mathrm{Re}(s) \geq 1/2$, while $L(s,\Delta_{\circ},BF)^{-1}$ has its zeros contained in the 
region $\mathrm{Re}(s) \leq 0$. Since half planes $\mathrm{Re}(s) \geq 1/2$ and $\mathrm{Re}(s) \leq 0$ are disjoint, they do not share factors in $\mathbb{C}[q^{\pm s}]$.
As $\rho$ is unitary, the poles of the product of $L$-factors in the numerator must lie on the line $\mathrm{Re}(s)=(\ell-i)/2$ for $i=0,\dotsm,\ell-2,\ell-1$, while those in the denominator 
are located on the line $\mathrm{Re}(s)=-i/2$ for $i=0,\dotsm,\ell-2,\ell-1,\ell$. Therefore they do not have common factors at all. We establish the desired identity from the observation that the ratios \eqref{reduce} and
$
  \gamma(s,\Delta_{\circ},BF,\psi) \sim \dfrac{L(1/2-s,\widetilde{\Delta}_{\circ},\delta_{-1/2},BF)}{L(s,\Delta_{\circ},BF)}
$
are all reduced and indices $i$'s are rearranged.
\end{proof}

Theorem \ref{BF-disc} is the key step to improve the factorization to the set of discrete series representations. 
If we can do this, then the application of the Langlands classification Theorem allows us to extend it 
to all irreducible admissible representations.

\begin{theorem}[Discrete series cases] 
\label{BF-disc}
Let $\Delta$ be an irreducible essentially square integrable representation of $GL_m$. Then
\[
L(s,\Delta,BF)=L(s+1/2,\Delta)L(2s,\Delta,\wedge^2).
\]
\end{theorem}

\begin{proof}
We choose an unramified quasi-character $\nu^{s_1}$, $s_1 \in \mathbb{C}$, so that $\Delta=\Delta_{\circ}\nu^{s_1}$, where $\Delta_{\circ}$ is an irreducible 
square integrable representation of $GL_m$. We can easily verify that $L(s,\Delta_{\circ}\nu^{s_1},BF)=L(s+s_1,\Delta_{\circ},BF)$, $L(s+1/2,\Delta_{\circ}\nu^{s_1})=L(s+s_1+1/2,\Delta_{\circ})$, and $L(2s,\Delta_{\circ}\nu^{s_1},\wedge^2)=L(2s+2s_1,\Delta_{\circ},\wedge^2)$. Hence for the calculation purpose,
we may assume that $\Delta=\Delta_{\circ}$ is unitary. The representation $\Delta$ is the segment consisting of supercuspidal representations of the form
$
 \Delta_{\circ}=[\rho_{\circ}\nu^{-\frac{\ell-1}{2}},\dotsm,\rho_{\circ}\nu^{\frac{\ell-1}{2}}],
$  
where $\rho_{\circ}$ is an irreducible unitary supercuspidal representation of $GL_r$ with $m=\ell r$.
We replace $\sigma$ by $\textbf{1}_{F^{\times}}$ in \cite[\S 2.6.2 Corollary]{CogPS}. Then the formula becomes
\[
 L(s+1/2,\Delta_{\circ})=L(s+1/2,\Delta_{\circ} \times \textbf{1}_{F^{\times}})=L( s+\ell/2,\rho_{\circ}).
\]
 This may also be seen from \cite[Proposition 3.1.3]{Jacquet}. Now we are left with invoking Theorem \ref{Discrete}.
\end{proof}

Finally, Theorem \ref{BF-all} renders the factorization result unconditional.

\begin{theorem}
\label{BF-all}
Let $\pi$ be an irreducible admissible representation of $GL_m$. Then we have 
\[
  L(s,\pi,BF)=L(s+1/2,\pi)L(2s,\pi,\wedge^2).
\]
\end{theorem}

\begin{proof}
We realize $\pi$ as the unique Langlands quotient of Langlands type $\Xi=\mathrm{Ind}^{GL_m}_{\rm Q}({\Delta_{\circ}}_1\nu^{u_1} \otimes {\Delta_{\circ}}_2\nu^{u_2} \otimes \dotsm \otimes {\Delta_{\circ}}_t\nu^{u_t})$ which is again of Whittaker type. The local Bump-Friedberg $L$-function is defined to be $ L(s,\pi,BF)=L(s,\Xi,BF)$. By \cite[Theorem 5.2]{Matringe15}, we have the equality
\[
  L(s,\Xi,BF)=\prod_{1 \leq k \leq t} L(s+u_k,{\Delta_{\circ}}_k,BF) \prod_{1 \leq i < j  \leq t} L(2s+u_i+u_j, {\Delta_{\circ}}_i \times {\Delta_{\circ}}_j).
\]
Applying Theorem \ref{BF-disc}, the product can be further decomposed as
\[
 L(s,\Xi,BF)=\prod_{1 \leq k \leq t} L(s+u_k+1/2,{\Delta_{\circ}}_k) 
 \prod_{1 \leq k \leq t} L(2s+2u_k,{\Delta_{\circ}}_k,\wedge^2) \prod_{1 \leq i < j  \leq t} L(2s+u_i+u_j, {\Delta_{\circ}}_i \times {\Delta_{\circ}}_j).
\]
Collecting the contributions for the first product $\prod L(s+u_k+1/2,{\Delta_{\circ}}_k) $ gives the standard $L$-factor $L(s+1/2,\Xi)=L(s+1/2,\pi)$ by \cite[Theorem 3.4]{Jacquet},
while gathering those for the rest product $\prod L(2s+2u_k,{\Delta_{\circ}}_k,\wedge^2) \prod L(2s+u_i+u_j, {\Delta_{\circ}}_i \times {\Delta_{\circ}}_j)$ yields the exterior square $L$-factor $L(2s,\Xi,\wedge^2)=L(2s,\pi,\wedge^2)$ by Theorem \ref{ext-langlands-prod}.
\end{proof}

We end this section with relating Bump-Friedberg $L$-factors to Galois theoretic counterparts. In conclusion, it is a consequence of the local Langlands correspondence that 
$L(s+1/2,\pi)=L( s+1/2, \phi(\pi))$ combined with Theorem \ref{exter-equal} and Theorem \ref{BF-all}.

\begin{theorem}
\label{BF-equal}
Let $\pi$ be an irreducible admissible representation of $GL_m(F)$ and $\phi(\pi)$ its associated Weil-Deligne representation of $W'_F$ under the local Langlands correspondence. Then we have
\[
  L(s,\pi,BF)=L( s+1/2, \phi(\pi) )L(2s,\wedge^2 ( \phi(\pi))).
\]
\end{theorem}

\subsection{Asai $L$-factors} 
\label{sec:Asai}

Let $E$ be a quadratic extension of $F$. We denote by $x \mapsto \overline{x}$ the non-trivial associated Galois action.
We fix an element $z \in E^{\times}$ such that $\overline{z}=-z$ and a non-trivial character $\psi_F$ of $F$. Let
\[
  \psi_E(x)=\psi_F \left( \frac{x-\overline{x}}{z-\overline{z}} \right), \quad x \in E.
\]
Then the additive character $\psi_E$ of $E$ is trivial on $F$ and defines a character of $N_m(E)$, which by abuse of notation we again denote by $\psi_E$. 
We shall use the Fourier transform induced by the additive character $\psi$ on the space of
Schwartz-Bruhat space $\mathcal{S}(F^m)$. Let $\pi=\mathrm{Ind}^{GL_m}_{\rm Q}(\Delta_1 \otimes \dotsm \otimes \Delta_t) \in \mathcal{A}_E(m)$ be a parabolically induced representation with an associated Whittaker model $\mathcal{W}(\pi,\psi_E)$. For each Whittaker function $W \in \mathcal{W}(\pi,\psi_E)$ and each Schwartz-Bruhat function
$\Phi \in \mathcal{S}(F^m)$, we define the local Flicker integral \cite{Fil88,Fil93} by
\[
  \mathcal{Z}(s,W,\Phi)=\int_{N_m \backslash GL_m} W(g)\Phi(e_mg) |\mathrm{det}(g)|^s dg
\]
which is absolutely convergent when the real part of $s$ is sufficiently large. Each $\mathcal{Z}(s,W,\Phi)$ is a rational function of $q^{-s}$
and hence extends meromorphically to all of $\mathbb{C}$. These integrals $\mathcal{Z}(s,W,\Phi)$  span a fractional ideal $\mathcal{I}(\pi,As)$.
of $\mathbb{C}[q^{\pm s}]$ generated by a normalized generator of the form $P(q^{-s})^{-1}$ where the polynomial $P(X) \in \mathbb{C}[X]$
satisfies $P(0)=1$. The {\it local Asai} $L$-{\it function} attached to $\pi$ is defined by such a unique normalized generator \cite[Definition 3.1]{Matringe09}:
\[
 L(s,\pi,As)=\frac{1}{P(q^{-s})}.
\]
Let us define the {\it local Asai} $\varepsilon$-factor by, as usual \cite[\S 3]{Matringe15} (Cf. \cite[\S 8]{AKMSS}):
\[
  \varepsilon(s,\pi,\psi,As)=\gamma(s,\pi,\psi,As) \frac{L(s,\pi,As)}{L(1-s,\pi^{\iota},As)}.
\]

\par
The Weil-Deligne group $W'_E$ of $E$ is of index two in the Weil-Deligne group $W'_F$ of $F$ and the quotient $W'_F \slash W'_E$ is naturally identified with $\mathrm{Gal}(E / F)$. We fix an element $\sigma$ in $W'_F$ which does not belong to $W'_E$ once for all. The image of $\sigma$ in  
$W'_F \slash W'_E$ is the non-trivial element of  $\mathrm{Gal}(E / F)$ which by abuse notation also denoted $\sigma$. Given a $m$-dimensional (complex) Frobenius semi-simple representation $\phi$ of $W'_E$, the {\it Asai representation} $As ( \phi) : W'_F \rightarrow {\rm GL}(\mathbb{C}^m \otimes \mathbb{C}^m) \simeq GL_{m^2}(\mathbb{C})$ given by (twisted) tensor induction of $\phi$ is defined as (Cf. \cite[\S 2.1]{AR}, \cite[\S 2]{Kri}, and \cite[\S 1.2]{Shankman}):
\[
  As(\phi)(\tau)(v \otimes w)=
  \begin{cases}
   \phi(\tau)(v) \otimes \phi(\sigma \tau \sigma^{-1} )(w)& \text{if $\tau \in W'_E$;}  \\
     \phi(\tau\sigma^{-1})(w) \otimes  \phi(\sigma\tau)(v)&  \text{if $\tau \notin W'_E$.}\\
  \end{cases}
\]
We then denote by $L(s,As( \phi))$ the {\it Artin} $L$-{\it factor} attached to Asai representation.

\par
The conjugation $\sigma$ extends naturally to an automorphism of $GL_m(E)$, which we also denote by $\sigma$. If $\pi \in \mathcal{A}_E(m)$, we denote by $\pi^{\sigma}$
the representation $g \mapsto \pi(\sigma(g))$.

\begin{theorem}
\label{equal-asai}
Let $\pi$ be an irreducible admissible representation of $GL_m(E)$ and $\phi(\pi)$ its associated Weil-Deligne representation of $W'_E$ under the local Langlands correspondence. Then we have
\[
  L(s,\pi,As)=\mathcal{L}(s,\pi,As)=L(s,As(\phi(\pi))).
\]
\end{theorem}
\begin{proof}
We first consider the case of irreducible unitary supercuspidal representations $\rho_{\circ}$ of $GL_r$. As a consequence of \cite{AKMSS} joined with \cite[Proposition 6]{AR} and \cite{HL-Asai}, we have
\[
 L(s,\rho_{\circ},As)=\mathcal{L}(s,\rho_{\circ},As).
\]

\par
Let $\Delta$ be a discrete series representation. In the spirit of twists by unramified characters for Langlands-Shahidi theoretic $L$-factors \cite[\S 3.1.(vi)]{HL-Asai} and Rankin-Selberg theoretic $L$-factors \cite[Theorem 2.3]{Matringe09}, there is no harm to assume that $\Delta=\Delta_{\circ}$ is an irreducible square integrable representation of $GL_{r\ell}$ associated to the segment $[\rho_{\circ}\nu^{-\frac{\ell-1}{2}},\dotsm,\rho_{\circ}\nu^{\frac{\ell-1}{2}}]$ with $\rho_{\circ}$ an irreducible unitary supercuspidal representation of $GL_r$. 
Let $\chi_{E \slash F}$ be an extension to $E^{\times}$ of the character $F^{\times}$ associated to $E \slash F$ by the local class field theory. As explained in \cite[Appendix A]{AKMSS},
\cite[Corollary 4.24]{Matringe09} and \cite[Theorem 4.26]{Matringe09} driven from the Cogdell and Piatetski-Shapiro method similar to \S \ref{deform-special} depend on the complete classification of $GL_m(F)$-distinguished representations \cite{Matringe11}. Looking at the proof of this Proposition, we need to check that the $GL_m(F)$-distinguished representation, namely, 
$\mathrm{Hom}_{GL_m(F)}(\pi,\textbf{1}_{GL_m(F)}) \neq \{ 0\}$, is still Galois self-dual, $\pi^{\sigma} \simeq \widetilde{\pi}$, for any non-archimedean local field of odd residual characteristic. It is 
presently written in this generality in the literature (Cf. \cite[\S 3.2.12]{Offen}). Counting on the weak multiplicativity of $\gamma(s,\pi,As,\psi)$ \cite[Corollary 4.24]{Matringe09}, we get the results below using arguments parallel to the one employed in the proof of Goldberg
\cite[Theorem 5.6]{Goldberg}:
\[
 L(s,\Delta_{\circ},As)=\prod_{i=1}^{\ell/2}L(s,\rho_{\circ}\nu^{(\ell+1)/2-i},As)L(s,\chi_{E \slash F} \otimes \rho_{\circ}\nu^{\ell/2-i},As)
\]
when $\ell$ is even and
\[
  L(s,\Delta_{\circ},As)=\prod_{i=1}^{(\ell+1)/2}L(s,\rho_{\circ}\nu^{(\ell+1)/2-i},As)\prod_{i=1}^{(\ell-1)/2}L(s,\chi_{E \slash F} \otimes \rho_{\circ}\nu^{\ell/2-i},As)
\]
when $\ell$ is odd. The expression is alike to that in Theorem \ref{Discrete}. This places us in a position to deduce
\begin{equation}
\label{asai-disc}
 L(s,\Delta,As)=\mathcal{L}(s,\Delta,As)
\end{equation}
for any discrete series representations $\Delta$ of $GL_{r\ell}$.

\par
In general, we realize $\pi$ as the unique Langlands quotient of Langlands type $\Xi=\mathrm{Ind}^{GL_m}_{\rm Q}({\Delta_{\circ}}_1\nu^{u_1} \otimes {\Delta_{\circ}}_2\nu^{u_2} \otimes \dotsm \otimes {\Delta_{\circ}}_t\nu^{u_t})$. As such, by the inductive relation of $L(s,\pi,As)$ \cite[Theorem 4.26]{Matringe09}, one has the equality
\[
 L(s,\pi,As)=\prod_{1 \leq k \leq t} L(s+2u_k,{\Delta_{\circ}}_k,\wedge^2) \prod_{1 \leq i < j  \leq t} L(s+u_i+u_j, {\Delta_{\circ}}_i \times {\Delta^{\sigma}_{\circ}}_j).
\]
Consequently, the equality
\[
  L(s,\pi,As)=\mathcal{L}(s,\pi,As)
\]
follows from \cite[\S 4.2]{HL-Asai} along with \eqref{asai-disc} for all irreducible admissible representations $\pi$ of $GL_m(E)$. The remaining part is simply to quote the main Theorem of Henniart and Lomel\'{\i} \cite[Theorem 4.3]{HL-Asai}.
\end{proof}


\begin{acknowledgements}
As this paper is a conclusion of my Ph.D. work, the author would like to thank, my advisor James Cogdell, for countless stimulus discussions and for proposing to
express $L$-factors in terms of those for supercuspidal representations. I sincerely thank to Nadir Matringe for patiently answering questions about linear and Shalika periods. I also thank to Shantanu Agarwal for drawing the author's attention to positive characteristic. In particular, I am indebted to Muthu Krishnamurthy for fruitful
mathematical communications over the years and describing a whole picture of the Langlands-Shahidi method to me.
I am grateful to the University of Iowa and the University of Maine for their hospitality and support while the article is written.
Lastly, I would like to thank to the anonymous referee for many valuable remarks and suggestions, which significantly improve exposition and organization of this paper.
\end{acknowledgements}

 \bibliographystyle{amsplain}

\begin{bibdiv}
\begin{biblist}

\bib{AKMSS}{article}{
   author={Anandavardhanan, U. K.},
   author={Kurinczuk, R.},
   author={Matringe, N.},
   author={S\'{e}cherre, V.},
   author={Stevens, S.},
   title={Galois self-dual cuspidal types and Asai local factors},
   journal={J. Eur. Math. Soc. (JEMS)},
   volume={23},
   date={2021},
   number={9},
   pages={3129--3191},
  }
  
  
  \bib{AR}{article}{
   author={Anandavardhanan, U. K.},
   author={Rajan, C. S.},
   title={Distinguished representations, base change, and reducibility for
   unitary groups},
   journal={Int. Math. Res. Not.},
   date={2005},
   number={14},
   pages={841--854},
}


\bib{BeZeSurvey}{article}{
   author={Bern\v{s}te\u{\i}n, I. N.},
   author={Zelevinski\u{\i}, A. V.},
   title={Representations of the group $GL(n,F),$ where $F$ is a local
   non-Archimedean field},
   language={Russian},
   journal={Uspehi Mat. Nauk},
   volume={31},
   date={1976},
   number={3(189)},
   pages={5--70},
}

\bib{BeZe}{article}{
   author={Bernstein, I. N.},
   author={Zelevinsky, A. V.},
   title={Induced representations of reductive ${\germ p}$-adic groups. I},
   journal={Ann. Sci. \'{E}cole Norm. Sup. (4)},
   volume={10},
   date={1977},
   number={4},
   pages={441--472},
}


\bib{BF}{article}{
   author={Bump, Daniel},
   author={Friedberg, Solomon},
   title={The exterior square automorphic $L$-functions on ${\rm GL}(n)$},
   conference={
      title={Festschrift in honor of I. I. Piatetski-Shapiro on the occasion
      of his sixtieth birthday, Part II},
      address={Ramat Aviv},
      date={1989},
   },
   book={
      series={Israel Math. Conf. Proc.},
      volume={3},
      publisher={Weizmann, Jerusalem},
   },
   date={1990},
   pages={47--65},
  }
  
  
  \bib{BD08}{article}{
   author={Blanc, Philippe},
   author={Delorme, Patrick},
   title={Vecteurs distributions $H$-invariants de repr\'{e}sentations induites,
   pour un espace sym\'{e}trique r\'{e}ductif $p$-adique $G/H$},
   language={French, with English and French summaries},
   journal={Ann. Inst. Fourier (Grenoble)},
   volume={58},
   date={2008},
   number={1},
   pages={213--261},
 }
  
  \bib{Chen-Gan}{article}{
  author={Chen, Rui},
  author={Gan, Wee Teck},
   title={\it Unitary Friedberg-Jacquet periods},
    pages={Preprint, available at \url{https://arxiv.org/abs/2108.04064}},
  }

\bib{Cogdell-DOC}{article}{
   author={Cogdell, J. W.},
   title={Analytic theory of $L$-functions for ${\rm GL}_n$},
   conference={
      title={An introduction to the Langlands program},
      address={Jerusalem},
      date={2001},
   },
   book={
      publisher={Birkh\"{a}user Boston, Boston, MA},
   },
   date={2003},
   pages={197--228},
 }


\bib{CoMa}{article}{
   author={Cogdell, James W.},
   author={Matringe, Nadir},
   title={The functional equation of the Jacquet-Shalika integral
   representation of the local exterior-square $L$-function},
   journal={Math. Res. Lett.},
   volume={22},
   date={2015},
   number={3},
   pages={697--717},
   }
   
   \bib{CP94}{article}{
   author={Cogdell, James W.},
   author={Piatetski-Shapiro, I. I.},
   title={Exterior square $L$-function for ${\rm GL}(n)$},
   journal={Talk given at the Fields Institute, April},
   date={1994},
    pages={available at \url{https://people.math.osu.edu/cogdell.1/exterior-www.pdf}},
      }

   
   \bib{CogPS}{article}{
   author={Cogdell, J. W.},
   author={Piatetski-Shapiro, I. I.},
   title={Derivatives and L-functions for $GL_n$},
   conference={
      title={Representation theory, number theory, and invariant theory},
   },
   book={
      series={Progr. Math.},
      volume={323},
      publisher={Birkh\"{a}user/Springer, Cham},
   },
   date={2017},
   pages={115--173},
  }
  
  \bib{Flath}{article}{
   author={Flath, D.},
   title={Decomposition of representations into tensor products},
   conference={
      title={Automorphic forms, representations and $L$-functions},
      address={Proc. Sympos. Pure Math., Oregon State Univ., Corvallis,
      Ore.},
      date={1977},
   },
   book={
      series={Proc. Sympos. Pure Math., XXXIII},
      publisher={Amer. Math. Soc., Providence, R.I.},
   },
   date={1979},
   pages={179--183},
}
 
 
 \bib{Fil88}{article}{
   author={Flicker, Yuval Z.},
   title={Twisted tensors and Euler products},
   language={English, with French summary},
   journal={Bull. Soc. Math. France},
   volume={116},
   date={1988},
   number={3},
   pages={295--313},
  }
 
  
  \bib{Fil93}{article}{
   author={Flicker, Yuval Z.},
   title={On zeroes of the twisted tensor $L$-function},
   journal={Math. Ann.},
   volume={297},
   date={1993},
   number={2},
   pages={199--219},
}

\bib{Gan}{article}{
   author={Gan, Wee Teck},
   title={Periods and theta correspondence},
   conference={
      title={Representations of reductive groups},
   },
   book={
      series={Proc. Sympos. Pure Math.},
      volume={101},
      publisher={Amer. Math. Soc., Providence, RI},
   },
   date={2019},
   pages={113--132},
}

\bib{Gan-L}{article}{
   author={Gan, Wee Teck},
   author={Lomel\'{\i}, Luis},
   title={Globalization of supercuspidal representations over function
   fields and applications},
   journal={J. Eur. Math. Soc. (JEMS)},
   volume={20},
   date={2018},
   number={11},
   pages={2813--2858},
   issn={1435-9855},
  }


\bib{GL}{article}{
   author={Ganapathy, Radhika},
   author={Lomel\'{\i}, Luis},
   title={On twisted exterior and symmetric square $\gamma$-factors},
   language={English, with English and French summaries},
   journal={Ann. Inst. Fourier (Grenoble)},
   volume={65},
   date={2015},
   number={3},
   pages={1105--1132},
  }
  
  \bib{GK}{article}{
   author={Gel\cprime fand, I. M.},
   author={Kajdan, D. A.},
   title={Representations of the group ${\rm GL}(n,K)$ where $K$ is a local
   field},
   conference={
      title={Lie groups and their representations},
      address={Proc. Summer School, Bolyai J\'{a}nos Math. Soc., Budapest},
      date={1971},
   },
   book={
      publisher={Halsted, New York},
   },
   date={1975},
   pages={95--118},
}

  
  \bib{Goldberg}{article}{
   author={Goldberg, David},
   title={Some results on reducibility for unitary groups and local Asai
   $L$-functions},
   journal={J. Reine Angew. Math.},
   volume={448},
   date={1994},
   pages={65--95},
 }
  
  
  \bib{Henniart}{article}{
   author={Henniart, Guy},
   title={Correspondance de Langlands et fonctions $L$ des carr\'{e}s ext\'{e}rieur
   et sym\'{e}trique},
   language={French},
   journal={Int. Math. Res. Not. IMRN},
   date={2010},
   number={4},
   pages={633--673},
  }
  
  \bib{HL-Exterior}{article}{
   author={Henniart, Guy},
   author={Lomel\'{\i}, Luis},
   title={Local-to-global extensions for ${\rm GL}_n$ in non-zero
   characteristic: a characterization of $\gamma_F(s,\pi,{\rm Sym}^2,\psi)$
   and $\gamma_F(s,\pi,\wedge^2,\psi)$},
   journal={Amer. J. Math.},
   volume={133},
   date={2011},
   number={1},
   pages={187--196},
  }
  
  \bib{HL-Asai}{article}{
   author={Henniart, Guy},
   author={Lomel\'{\i}, Luis},
   title={Characterization of $\gamma$-factors: the Asai case},
   journal={Int. Math. Res. Not. IMRN},
   date={2013},
   number={17},
   pages={4085--4099},
   }
  
  \bib{HL-RS}{article}{
   author={Henniart, Guy},
   author={Lomel\'{\i}, Luis},
   title={Uniqueness of Rankin-Selberg products},
   journal={J. Number Theory},
   volume={133},
   date={2013},
   number={12},
   pages={4024--4035},
  }
    
  \bib{Jacquet}{article}{
   author={Jacquet, Herv\'{e}},
   title={Principal $L$-functions of the linear group},
   conference={
      title={Automorphic forms, representations and $L$-functions},
      address={Proc. Sympos. Pure Math., Oregon State Univ., Corvallis,
      Ore.},
      date={1977},
   },
   book={
      series={Proc. Sympos. Pure Math., XXXIII},
      publisher={Amer. Math. Soc., Providence, R.I.},
   },
   date={1979},
   pages={63--86},
}
  
  \bib{JPSS}{article}{
   author={Jacquet, H.},
   author={Piatetskii-Shapiro, I. I.},
   author={Shalika, J. A.},
   title={Rankin-Selberg convolutions},
   journal={Amer. J. Math.},
   volume={105},
   date={1983},
   number={2},
   pages={367--464},
  }
  
  \bib{JaSh81}{article}{
   author={Jacquet, H.},
   author={Shalika, J. A.},
   title={On Euler products and the classification of automorphic
   representations. I},
   journal={Amer. J. Math.},
   volume={103},
   date={1981},
   number={3},
   pages={499--558},
   }
      
   \bib{JaSh88}{article}{
   author={Jacquet, Herv\'{e}},
   author={Shalika, Joseph},
   title={Exterior square $L$-functions},
   conference={
      title={Automorphic forms, Shimura varieties, and $L$-functions, Vol.
      II},
      address={Ann Arbor, MI},
      date={1988},
   },
   book={
      series={Perspect. Math.},
      volume={11},
      publisher={Academic Press, Boston, MA},
   },
   date={1990},
   pages={143--226},
 }
 
 \bib{Jo20}{article}{
   author={Jo, Yeongseong},
   title={Factorization of the local exterior square $L$-function of $GL_m$},
   journal={Manuscripta Math.},
   volume={162},
   date={2020},
   number={3-4},
   pages={493--536},
  }
  
  \bib{Jo20-2}{article}{
   author={Jo, Yeongseong},
   title={Derivatives and exceptional poles of the local exterior square
   $L$-function for $GL_m$},
   journal={Math. Z.},
   volume={294},
   date={2020},
   number={3-4},
   pages={1687--1725},
}

\bib{Jo21}{article}{
   author={Jo, Yeongseong},
   title={Rankin-Selberg integrals for local symmetric square factors on
   $GL(2)$},
   journal={Mathematika},
   volume={67},
   date={2021},
   number={2},
   pages={388--421},
  }


\bib{Kable}{article}{
   author={Kable, Anthony C.},
   title={Asai $L$-functions and Jacquet's conjecture},
   journal={Amer. J. Math.},
   volume={126},
   date={2004},
   number={4},
   pages={789--820},
  }
  
  \bib{Kaplan}{article}{
   author={Kaplan, Eyal},
   title={The characterization of theta-distinguished representations of
   ${\rm GL}(n)$},
   journal={Israel J. Math.},
   volume={222},
   date={2017},
   number={2},
   pages={551--598},
}

\bib{Kri}{article}{
   author={Krishnamurthy, M.},
   title={Determination of cusp forms on $GL(2)$ by coefficients restricted
   to quadratic subfields (with an appendix by Dipendra Prasad and Dinakar
   Ramakrishnan)},
   journal={J. Number Theory},
   volume={132},
   date={2012},
   number={6},
   pages={1359--1384},
}

\bib{KeRa}{article}{
   author={Kewat, Pramod Kumar},
   author={Raghunathan, Ravi},
   title={On the local and global exterior square $L$-functions of ${\rm
   GL}_n$},
   journal={Math. Res. Lett.},
   volume={19},
   date={2012},
   number={4},
   pages={785--804},
  }
  
  \bib{Lomeli}{article}{
   author={Lomel\'{\i}, Luis Alberto},
   title={On automorphic $L$-functions in positive characteristic},
   journal={Ann. Inst. Fourier (Grenoble)},
   volume={66},
   date={2016},
   number={5},
   pages={1733--1771},
   }
  
 \bib{Matringe09}{article}{
   author={Matringe, Nadir},
   title={Conjectures about distinction and local Asai $L$-functions},
   journal={Int. Math. Res. Not. IMRN},
   date={2009},
   number={9},
   pages={1699--1741},
   }
   
   \bib{Matringe11}{article}{
   author={Matringe, Nadir},
   title={Distinguished generic representations of ${\rm GL}(n)$ over
   $p$-adic fields},
   journal={Int. Math. Res. Not. IMRN},
   date={2011},
   number={1},
   pages={74--95},
   }
 
   \bib{Matringe}{article}{
   author={Matringe, Nadir},
   title={Linear and Shalika local periods for the mirabolic group, and some
   consequences},
   journal={J. Number Theory},
   volume={138},
   date={2014},
   pages={1--19},
 }
 
 \bib{Matringe15}{article}{
   author={Matringe, Nadir},
   title={On the local Bump-Friedberg $L$-function},
   journal={J. Reine Angew. Math.},
   volume={709},
   date={2015},
   pages={119--170},
}

\bib{Matringe17}{article}{
   author={Matringe, Nadir},
   title={Shalika periods and parabolic induction for $GL(n)$ over a
   non-archimedean local field},
   journal={Bull. Lond. Math. Soc.},
   volume={49},
   date={2017},
   number={3},
   pages={417--427},
   }
   
   \bib{Offen}{article}{
   author={Offen, Omer},
   title={Period integrals of automorphic forms and local distinction},
   conference={
      title={Relative aspects in representation theory, Langlands
      functoriality and automorphic forms},
   },
   book={
      series={Lecture Notes in Math.},
      volume={2221},
      publisher={Springer, Cham},
   },
   date={2018},
   pages={159--195},
  }


\bib{Rodier}{article}{
   author={Rodier, Fran\c{c}ois},
   title={Whittaker models for admissible representations of reductive
   $p$-adic split groups},
   conference={
      title={Harmonic analysis on homogeneous spaces},
      address={Proc. Sympos. Pure Math., Vol. XXVI, Williams Coll.,
      Williamstown, Mass.},
      date={1972},
   },
   book={
      publisher={Amer. Math. Soc., Providence, R.I.},
   },
   date={1973},
   pages={425--430},
  }
   
    \bib{Shankman}{article}{
  author={Shankman, Daniel}, 
   title={\it Local Langlands correspondence for Asai L
functions and $\epsilon$-factors},
    pages={Preprint, available at \url{https://arxiv.org/abs/1810.11852}},
  }

   
   \bib{Shahidi}{article}{
   author={Shahidi, Freydoon},
   title={Twisted endoscopy and reducibility of induced representations for
   $p$-adic groups},
   journal={Duke Math. J.},
   volume={66},
   date={1992},
   number={1},
   pages={1--41},
   }
   
   \bib{Tate}{article}{
   author={Tate, J.},
   title={Number theoretic background},
   conference={
      title={Automorphic forms, representations and $L$-functions},
      address={Proc. Sympos. Pure Math., Oregon State Univ., Corvallis,
      Ore.},
      date={1977},
   },
   book={
      series={Proc. Sympos. Pure Math., XXXIII},
      publisher={Amer. Math. Soc., Providence, R.I.},
   },
   date={1979},
   pages={3--26},
  }


\bib{Yamana}{article}{
   author={Yamana, Shunsuke},
   title={Local symmetric square $L$-factors of representations of general
   linear groups},
   journal={Pacific J. Math.},
   volume={286},
   date={2017},
   number={1},
   pages={215--256},
}

\bib{Yan22}{article}{
   author={Yang, Chang},
   title={Linear periods for unitary representations},
   journal={Math. Z.},
   volume={302},
   date={2022},
   number={4},
   pages={2253--2284},
  }
 
 \bib{Zelevinsky}{article}{
   author={Zelevinsky, A. V.},
   title={Induced representations of reductive ${\germ p}$-adic groups. II.
   On irreducible representations of ${\rm GL}(n)$},
   journal={Ann. Sci. \'{E}cole Norm. Sup. (4)},
   volume={13},
   date={1980},
   number={2},
   pages={165--210},
  }   
 
 
  

\end{biblist}
\end{bibdiv}

\end{document}